\pdfoutput=1

\documentclass[11pt]{article}
\usepackage{amssymb}
\usepackage{graphicx}
\parskip \baselineskip

\newtheorem{lm}{Lemma}

\newtheorem{prop}{Proposition}

\newtheorem{corol}{Corollary}

\newtheorem{ex}{Example}

\newtheorem{re}{Remark}

\newcommand{\ds}{\displaystyle}

\newcommand{\N}{\mathbb{N}}

\newcommand{\CL}{\cal L}

\newcommand{\ol}{\overline}

\newcommand{\ra}{\rightarrow}

\date{}

\begin{document}

\title{Tight embedding of modular lattices into partition lattices: progress and program}

\author{Marcel Wild }

\maketitle

\begin{quote}
A{\scriptsize BSTRACT}: {\footnotesize A famous Theorem of Pudlak and T{\' u}ma states that each finite lattice $L$ occurs as sublattice of a finite partition lattice. Here we derive, for modular lattices $L$, necessary and sufficient conditions for cover-preserving embeddability. Aspects of our work relate to Bjarni J{\' o}nsson.}
\end{quote}

\section{Introduction}

 A bird's eye view of this lengthy article follows in 1A. The detailed Section break-up comes in 1B. How Bjarni J{\' o}nsson's work relates to our article, is outlined at the end of the Introduction in 1C. 

{\bf 1A Bird's eye view.} Without further mention, all structures are assumed to be {\it finite}.
Every concept not explained in this article  is standard and can e.g. be explored in  [2] (lattices) or [10] (matroids).  
By definition a {\it tight} embedding $f: L\ra L'$
between lattices is a cover-preserving lattice homomorphism with $f(0)=0$. It is easy to see that tight implies injective. What is more, when both $L$ and $L'$ are graded then the height $d(L)$ equals $d(Im(f))$. We like to keep the definition of 'tight' loose enough to allow for $d(Im(f))<d(L')$. Let $Part(n)$ be the semimodular (whence graded) lattice of all  partitions of the set $[9]:=\{1,2,\ldots,n\}$. As is well known, when $L$ wants to tightly embed into $Part(n)$, then $L$ must be semimodular itself.  

About twenty five years ago I made strides towards finding necessary and sufficient conditions for a {\it modular} lattice $L$ to be tightly embeddable  into  $Part(n)$. This is because modular lattices enjoy a much richer structure theory than merely semimodular ones. One modular lattice $L_0$ together with a tight embedding into $Part(5)$ is given in Figure 1a (= 1(a) ).
For instance, 13, 25, 4  is shorthand for the partition $\{\{1,3\},\{2,5\},\{4\}\}$. 

\includegraphics[scale=0.43]{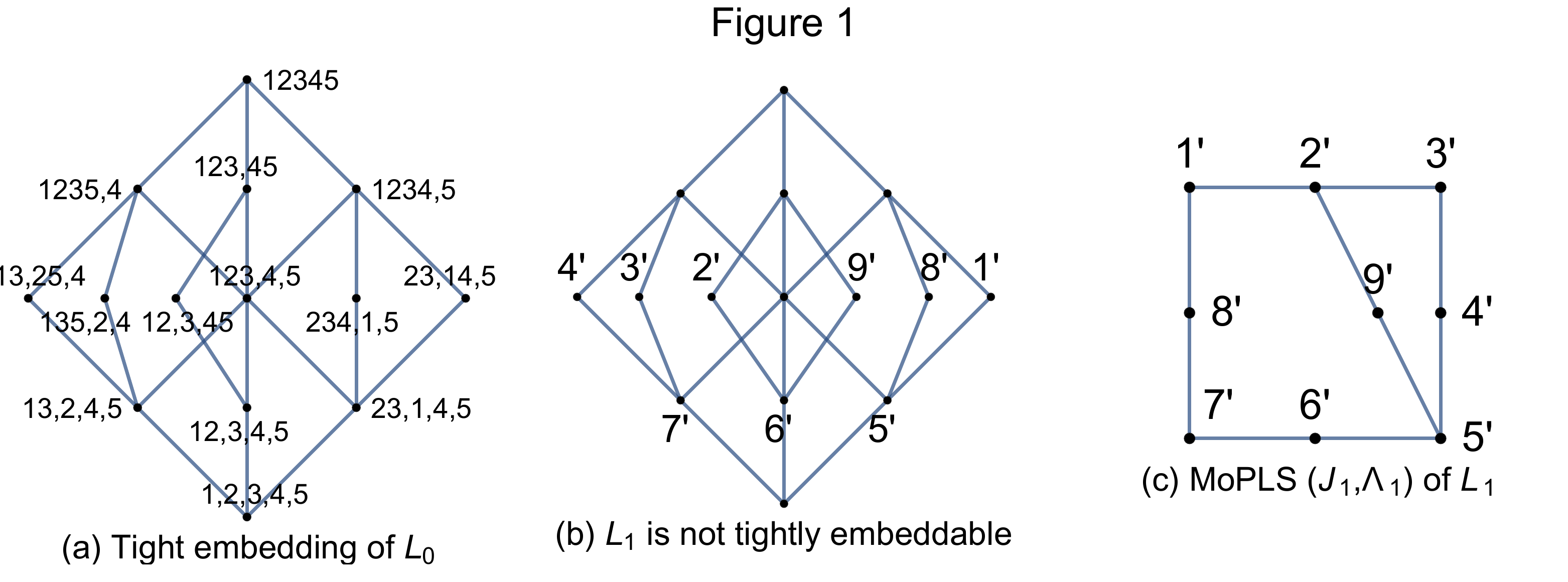}

Let us indicate why the one element larger lattice $L_1$ in Figure 1b does not admit a tight embedding into $Part(5)$. First off, there is no relation between the numbers $1,\ldots, 5$ and $1',\ldots,9'$ in Figures 1a and 1b. One checks that up to obvious symmetry Figure 1a yields the {\it only} tight embedding of $L_0$. Now the only partitions in $Part(5)$ lying between $\{12,3,4,5\}$ and $\{123,45\}$ are $\{123,4,5\}$ and $\{12,3,45\}$. Hence the element $9'\in L_1$ in Figure 1b cannot be assigned to a partition in $Part(5)$. What about $Part(n)$ for $n>5$? This does not help either, as implied by the results in Subsection 5B1.

Generally the necessary conditions for tight  embeddability into $Part(n)$ turned out [15] to be quite close to the sufficient ones.
To close the gap roughly speaking half of the remaining work (i.e. strengthening the necessary conditions) is lattice-theoretic, whereas the other half (i.e.  softening the sufficient conditions) is purely combinatorial. It is the latter half where the present article makes  progress.
This is mirrored by the fact that only three out of eight open questions (Questions 3,7,8)  involve modular lattices at all. In contrast Questions 2,4,5,6 are about graphs and binary matroids and e.g. ask whether certain types of binary matroids must in fact be graphic. We thus keep on soliciting the input of combinatorists; the contributions so far of Manoel Lemos and Jim Geelen are gratefully acknowledged and will be discussed in 8B, 8C. 

Progress initially sprang from the author's effort to present the key concepts and results of [15] in crisper ways. All of this harks back to a result of Huhn [6] that characterizes 2-distributive modular lattices in terms of certain forbidden sublattices. Huhn's paper sparked the groundbreaking paper of J{\'o}nsson and Nation [7] which showed that for each field $k$  each 2-distributive modular lattice $L$ with $|L|<|k|$ tightly embeds into the subspace lattice $L(k^n)\ (n=d(L))$. In turn [7] triggered [5], which triggered [15] and whence the present article.

 {\bf 1B The Section break-up.} In Section 2 we point out two straightforward sufficient conditions for tight partition embeddability of a modular lattice $L$. One is the tight embeddability of the subdirectly irreducible factors of $L$. The other condition (a consequence of the first, but easily shown directly) is the distributivity of $L$. While distributivity (=1-distributivity) is sufficient, the lesser known and much softer constraint of 2-distributivity turns out to be a necessary condition.   The next two Sections omit modular lattices altogether.

Section 3 reviews the definition of a partial linear space (PLS) and familiar concepts such as connected components and cycles. In fact only PLSes all of whose lines have cardinality 3 will be considered in our article. For instance the PLS in Figure 1c has point set $\{1',2',...,9'\}$ and five such lines. The 'square' that stands out is a cycle $C$ in the PLS-sense of the word. There is a path $P$ (in fact a single 
line $\{2',9',5'\}$) that connects the midpoint $2'$ of the $C$-line $\{1',2',3'\}$ with the $C$-junction $5'$. We shall strive to avoid such 
'non-benign' paths $P$. In contrast, useful Lemmata can be proven about BMPLes, i.e. PLSes that feature only 'benign' midpoint links.

 Section 4 links matroids [10] to partial linear spaces. 
Most of the matroid-related content in Section 4 and onwards actually is about graphs (i.e. graphic matroids), the rest about binary matroids and arbitrary matroids. We say that a simple matroid $M=M(E)$ with universe $E$ {\it line-preserving models} a PLS $(J,\Lambda)$ with point set $J$ and line-set $\Lambda$ if there is a bijection $\psi: J\ra E$ such that $\psi(\ell)$ is dependent in $M$ for each $\ell\in\Lambda$ (dependency condition of the first kind). If the matroid-rank $mrk(E)$ coincides with the PLS-rank  $rk(J,\Lambda)=|J|-|\Lambda|$ (rank condition of the first kind) then $\psi$ is a {\it rank-modeling} of $(J,\Lambda)$. Some
intriguing interplay arises between cycles in PLSes on the one hand, and circuits in graphs on the other. How Sections 3 and 4 relate to modular lattices only dawns in Section 5, and becomes fully apparent in Section 6. 

Lemma 5.1 in Section 5 states the following. Let $L$ be a modular lattice of length $d(L)$ and with set $J$ of join irreducibles. (For instance $L_1$ in Figure 1b has $J(L_1)=\{1',...,9'\}$.) Put $J(a):=\{p\in J:\ p\le a\}$ for all $a\in L$. Let $M(K)$ be a matroid with a lattice $LM(K)$ of closed subsets of $K$. Then there is a tight embedding $\Phi:L\ra LM(K)$ iff there is an injection $\varphi: J\ra K$ such that the induced submatroid $M(\varphi(J))$ of $M(K)$ is simple and there is an injection $\varphi$ that {\it lattice-models} $L$ in the following  sense. All sets $\varphi(J(a))\ (a\in L)$ must be closed in $M(\varphi(J))$ (dependency condition of the second kind). Furthermore, $mrk(\varphi(J))=d(L)$ (rank condition of the second kind). Specializing the arbitrary matroids in Lemma 5.1 to  graphic matroids  gives rise to Theorem 5.4 which states that a modular lattice $L$ with $d(L)=n$ is tightly embeddable into $Part(n+1)$ iff there is a graph on $n+1$ vertices that lattice-models $L$. This is all good and well, except that more verifiable conditions for tight embeddability are desirable.
The theory of 'abstract' PLSes developed in Sections 3 and 4 will go a long way towards fulfilling that wish. 

Namely, Section 6 puts the pieces together as follows. By our previous efforts we know how to construct  graphs $G$ modeling  PLSes $(J,\Lambda)$ in various ways. But since 
Theorem 5.4 demands $G$ to rather  {\it lattice-model}
  $L$, our $(J,\Lambda)$ better be linked to $L$ somehow. In a nutshell, Section 6 achieves this by forcing $J=J(L)$.
	At this point some readers may think about the Fundamental Theorem  of Projective Geometry. It links complemented modular lattices $L$ with particular types of PLSes,
i.e. projective spaces $(J,\Lambda)$. Here $J$ can be identified with the set of height 1 elements (=atoms) of $L$, and the set of lines $\Lambda$ is in bijection with the height 2 elements of $L$.
However a complemented modular lattice $L$ is not (except for trivial cases)  tightly partition embeddable, and so classic projective spaces won't do. Fortunately, if complementation lacks, projective spaces can be generalized to  PLSes on the {\it posets} $(J,\le)$  of join irreducibles. This has been pioneered by Herrmann and his co-authors [1],[5],[4]. We shall use the term MoPLS for a PLS that originates in this way from a modular lattice. One MoPLS of $L_1$ is shown in Figure 1c. Observe the crucial fact that any two elements on a line yield the same join. In contrast to the classic case usually plenty of MoPLSes can describe any fixed modular lattice $L$.  

Most important for us are {\it thin} modular lattices, i.e. 2-distributive and without covering sublattice $M_4$. For them the satisfaction of the rank condition of the first kind amounts to the satisfaction of the rank condition of the second kind. Unfortunately, it remains an open question (Question 3) whether  the  dependency condition of the first kind implies the dependency condition of the second kind. It was known [15] that the carry-over works under the additional proviso that the MoPLS is modeled  in a 'circuit-friendly' manner. In the present article (Theorem 6.3) we extend the class of MoPLSes for which circuit-friendly modeling can be guaranteed. The sufficiency of circuit-friendliness is based on the lengthy  lattice-theoretic Lemma 20 in [15]. This Lemma likely has a shorter proof, and perhaps its scope can be widened as well (Question 7).

Whatever the answer to Question 7, in Subsection 6C we launch an altogether different attack that shifts the burdon from 'modular latticians' to matroid theoreticians. The plan boils down to show that certain binary matroids  are actually graphic. Specifically, if $(J,\Lambda)$ is any abstract PLS with the property that each binary matroid that line-preserving rank-models $(J,\Lambda)$ must be graphic, then we call $(J,\Lambda)$ a {\it graph-trigger}. 
As  proof of concept the author shows in Theorem 6.6 that acyclic PLSes are graph-triggers. It is hoped that researchers more knowledgeable in matroid theory can handle more relevant types of PLSes. After all, finding conditions under which binary matroids (admittedly, not entangled with PLSes) must be graphic has a long history, starting with Tutte 1959 [13]. 

In order not to interrupt the storyline too much the longer proofs of Lemmas and Theorems will  be postponed to Section 7. This entails that the numbering of Figures in the main text may undergo apparent gaps. Section 8 houses further remarks that tend to be technical but are deemed  useful for future research. Throughout the article 'iff' means 'if and only if'. Some definitions (usually the more fundamental ones like {\bf Definition 2.2}) are announced in boldface, but others (like {\it path}) are simply set out in Italics within the plain text. Announced definitions, and also Proofs and Examples, always end with $\Box$.

{\bf 1C Connections to Bjarni J{\' o}nsson.} The impact of the  J{\' o}nsson-Nation paper [7] was mentioned already. We note in passing that [16], which briefly comes up in 2C, was invited by J{\' o}nsson. In the sequel of 1C (and only there) lattices need not be finite. In 1946 Whitman proved that every lattice $L$ has an embedding $f:L\ra Part(A)$. In 1953 J{\' o}nsson showed that this embedding can be chosen to be of {\it Type 3}. Roughly speaking this concerns the complexity, for all $x,\ y$ in $L$, of the equivalence relation (=partition) $f(x\vee y)$ in relation to $f(x)$ and $f(y)$. J{\' o}nsson further proved that a certain {\it Type 2} embedding is possible iff $L$ is modular. What is more, J{\' o}nsson showed that {\it Type 1} (i.e. permuting equivalence relations) defines a strict subclass of modular lattices, the so-called Arguesian lattices. See chapter 4 in [9] for more details on the Type 1,2,3 story. The next 'craze', i.e. whether the groundset $A$ in $Part(A)$ can be chosen finite when $L$ is finite, was settled affirmatively by Pudlak and T{\' u}ma in 1980. Unfortunately $|A|$ is super-exponential with respect to $|L|$, i.e. the embedding is everything but tight. The present article attempts to cure this state of affairs for the natural class of thin  lattices $L$.

\section{First steps}

Apart from the survey rendered in 1B we mention that in 2C we state a numerical inequality that holds in all modular lattices. Its sharpness is sufficient for tight partition embeddability. 

					{\bf 2A Distributive lattices.} An {\it order ideal} in a poset $(P,\le)$ is a subset $X\subseteq P$ such that from $y\le x\in X$ follows $y\in X$. The set $D(P,\le)$ of all order ideals is closed under $\cap$ and $\cup$, whence it is a (necessarily distributive) sublattice of the powerset lattice $\mathcal{P}(P)$. We denote by $J(L)$ the set of nonzero join-irreducibles of a lattice $L$. For each $a\in L$ put $J(a):=\{p\in J(L):\ p\le a\}$, and write $(J,\le)$ for the poset  arising from restricting the lattice ordering to the subset $J=J(L)$. By Birkhoff's Theorem each distributive lattice $L$ is isomorphic to $D(J,\le)$ via $a\mapsto J(a)$. Since therefore $L$ is isomorphic to  sublattice of $\mathcal{P}(J)$, it readily follows that $L$ tightly embeds into $Part(|J|)$. Namely, for $A\subseteq J$ put $f(A):=\{\{x\}:\ x\in J\setminus A\}\cup\{A\}$. Then $A\mapsto f(A)$ is a tight lattice embedding $\mathcal{P}(J)\ra Part(|J|)$.

{\bf 2B Subdirect products.} Let $s=s(N)$ be the number of maximal congruences $\theta_i$ of a modular lattice $N$. Letting $N_i:=N/\theta_i$  it is well known that $d(N)=d(N_1)+\cdots +d(N_s)$ and that $L$ tightly embeds into $N_1\times\cdots\times N_s$ as a subdirect product. Theorem 2.1 below is thus plausible but for the sake of completeness a full proof is provided in Section 7.

{\bf Theorem 2.1:} If the subdirectly irreducible factors $N_i$ of the modular lattice $N$ are tightly partition embeddable, then so is $N$. 

Notice that Theorem 2.1 covers the tight embeddability of distributive lattices since they are subdirect products of 2-element lattices $D_2$. As usual  denote by $M_n\ (n\ge 3)$ the length two modular lattice with $n$ join irreducibles. Since $M_3$ is tightly partition embeddable (see Figure 2d), each subdirect product $N$ of lattices $D_2$ and $M_3$ is as well.  For $N\subseteq M_3\times D_2$ the details are illustrated in Figure 2. 

If the modular lattice $FM(P,\le)$ freely generated by the poset $(P,\le)$ happens to be finite, then $FM(P,\le)$ is  a subdirect product of some $m$ lattices $M_3$ and $d$ lattices $D_2$  [17]. Hence $FM(P,\le)$ tightly embeds into $Part(n+1)$ for $n:=2m+d$.

\includegraphics[scale=0.3]{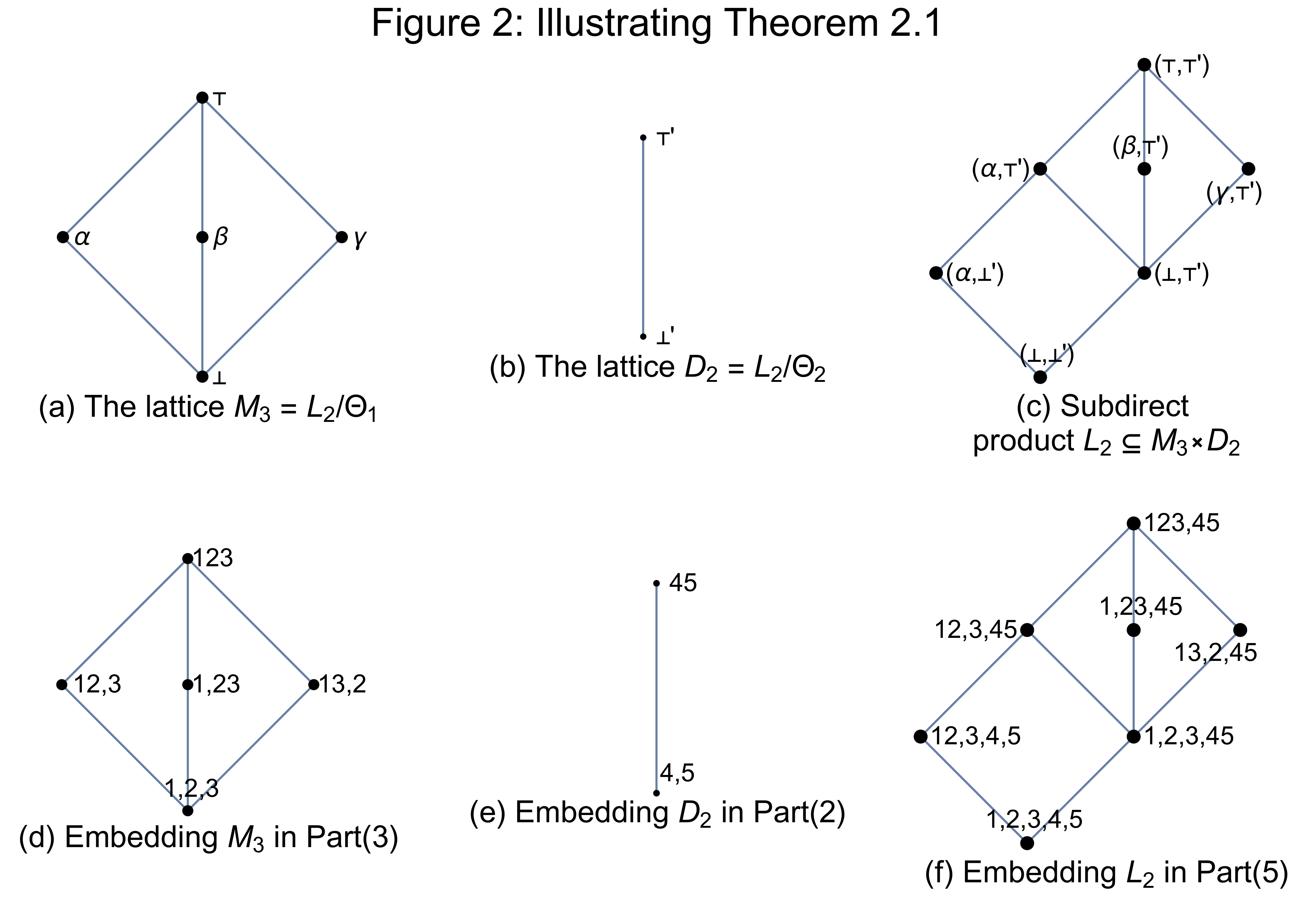}

{\bf 2C An arithmetic sufficient condition.} One can show [16] that each modular lattice $N$ satisfies 

(2.1)\qquad $|J(N)|\ge 2d(N)-s(N)$.

 For instance for $N=L_1$ in Figure 1b this becomes $9\ge 2\cdot 4-1$. For  distributive lattices $N$ inequality (2.1) is sharp because in fact $|J(N)|=d(N)=s(N)$. Since sharpness carries over to subdirect products, and sharpness holds for $M_3$ (i.e. $3=2\cdot 2 -1$), it likewise follows that (2.1) is sharp for $FM(P,\le)$  from above.  A complete characterization of the modular lattices for which (2.1) is sharp follows in 6A. They all are tightly partition embeddable.

{\bf 2D Thin lattices.}  It is well known that each interval of $Part(n)$ is isomorphic to a direct product of partition lattices. In particular a length 2 interval is isomorphic to $Part(3)=M_3$ or to $Part(2)\times Part(2)=D_2\times D_2$. Thus $M_4$ cannot be a covering sublattice of $Part(n)$. A length 3 interval of $Part(n)$ is isomorphic to either $Part(4)$ or $Part(3)\times
Part(2)$ or $Part(2)\times Part(2)\times Part(2)$. Since none of these lattices has more than 6 atoms, none of them is isomorphic to the subspace lattice of a 
nondegenerate projective plane since the smallest such plane is the Fano-plane in Figure 11a (Section 4) with 7 points. Its subspace lattice is isomorphic to $L(k^3)$, where $k=GF(2)$ is the 2-element field.

{\bf Definition 2.2.} A lattice $L$ is called {\it 2-distributive} if

(2.2)\qquad $a\wedge(b\vee c\vee d)=(a\wedge(b\vee c))\vee (a\wedge(b\vee d))\vee (a\wedge(c\vee d))$

and the dual identity hold for all $a,b,c,d\in L$.$\Box$

 When $L$ is modular, which is always the case for us, then (2.2) implies the dual identity. Evidently (2.2) is a (sweeping) generalization of the distributive law. As shown by Huhn   $L$ is 2-distributive iff it does not contain a length 3 interval isomorphic to the subspace lattice of a 
nondegenerate projective plane. It thus follows that a modular lattice that tightly embeds into a partition lattice must be {\it thin} in the sense of being 2-distributive and avoiding covering sublattices $M_4$.  Nevertheless, Figure 1b shows that some thin lattices stubbornly refuse to be tightly embedded into a partition lattice.

\section{Partial linear spaces on their own}

{\bf Definition 3.1:} A {\it partial linear space (PLS)} is an ordered pair $(J,\Lambda)$ consisting of a set $J$ of points and a set $\Lambda$ of 3-element subsets $l\subseteq \Lambda$ called {\it lines} such that 

(3.1)\quad $|l\cap l'| \le 1$ for all distinct $l,\  l'\in\Lambda$. $\Box$

Usually also lines of cardinality $>3$ are considered, but for us only cardinality 3 matters. As will be seen, this relates to the fact that $M_3$ is tightly embeddable into $Part(n)$ but $M_4$ is not.
In 3A we introduce paths and cycles in PLSes in unsurprising ways. In 3B to 3D we present increasingly complex types of PLSes, i.e. QIMPes, UMPes, and BMPLes. Sections 3E and 3F are  about iteratively piling up lines and paths, and 3G introduces the rank of a PLS smoothly as $rk(J,\Lambda):=|J|-|\Lambda|$.

{\bf 3A Paths and cycles.} Because of condition (3.1) any distinct points $p,\ q$ of a PLS $(J,\Lambda)$ lie on at most one common line which we then denote by $[p,q]$. For $n\ge 2$ a tuplet $P=[p_1,p_2,..,p_n]$ is called a {\it path} in $(J,\Lambda)$ if all lines $[p_i,p_{i+1}]\ (1\le i< n)$ exist, are distinct,
and for all $1\le i<j<n$ it holds that $[p_i,p_{i+1}]\cap[p_j,p_{j+1}]\neq \emptyset$ iff $j=i+1$. Consequently
the underlying point set $P^{*}:=[p_1,p_2]\cup\cdots\cup [p_{n-1},p_n]$ has cardinality $2n-1$.
Two points $p,q$ are {\it connected} if $p=q$ or there is a path $[p,...,q]$. This yields an equivalence relation whose $c(J,\Lambda)$ many classes are the {\it connected components} of $(J,\Lambda)$. Call $p$ {\it isolated} when $\{p\}$ is a connected component. The PLS in Figure 3b has six connected components, four of which are isolated points.

Albeit 'cycle' is used in graph theory, for us 'cycle' always refers to PLSes. (We shall soon be concerned with the corresponding structure in graphs, which  we name 'circuits'). Namely, a {\it cycle} C=$(p_1,p_2,..,p_n)$  is a path $[p_1,p_2,..,p_n]$ such that the line $[p_n,p_1]$ exists and features a new point $q$, i.e.
$C^{*}:=[p_1,\cdots,p_n]^*\cup\{q\}$ has cardinality $2n$.
Thus $(1,2,5)$ is a cycle in the PLS depicted in Figure 3a. In contrast $[1,2,3,4]$ is a path which is {\it no} cycle; 
albeit $[4,1]=\{4,1,5\}$ exists, the point $5$ is {\it not new} since it belongs to a preceeding line, i.e. $[2,3]$.

\includegraphics[scale=0.55]{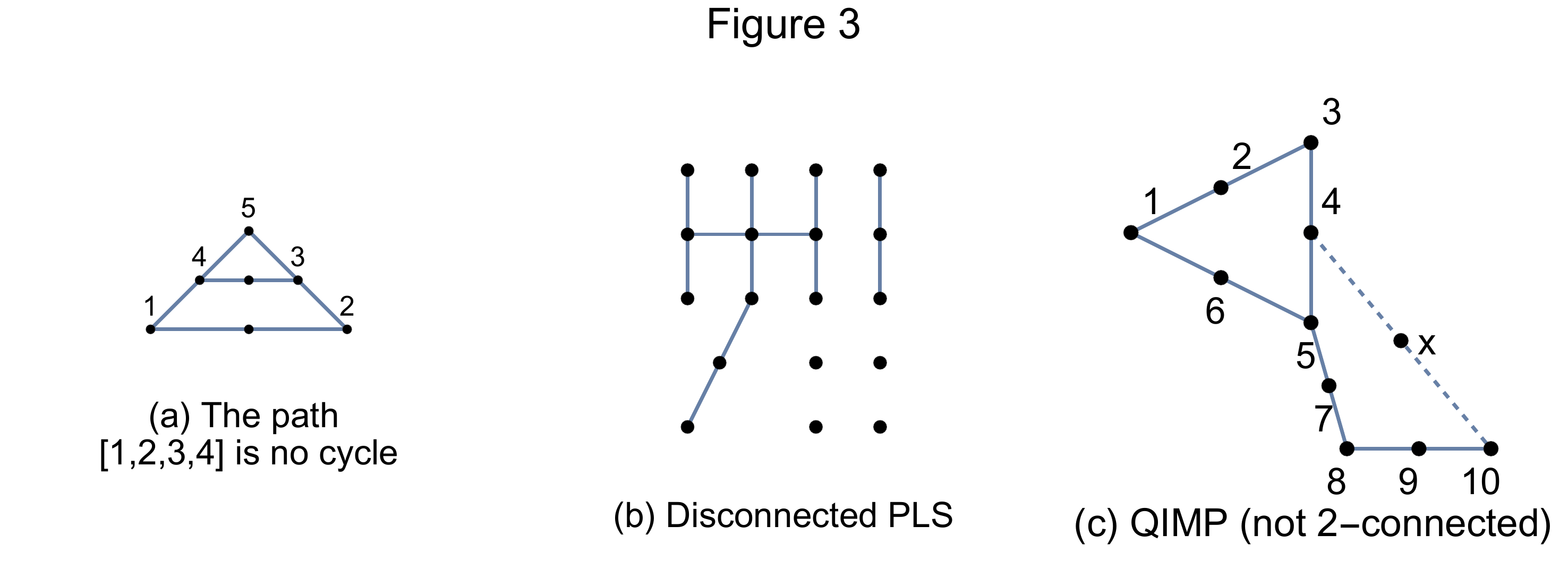}

 Necessarily each cycle $C=(p_1,p_2,..,p_n)$ has $n\ge 3$ and we call the points $p_i$ the {\it C-junctions}. For each {\it C-line}, i.e. $l=[p_i,p_{i+1}]$ (where $n+1:=1$), the unique point in $l\setminus[p_i,p_{i+1}]$ is called the {\it C-midpoint} of $\ell$. Thus $C$-lines can only intersect in a $C$-junction, not in a $C$-midpoint. A line is {\it interesting} if it is a $C$-line for at least one cycle $C$. A PLS without cycles is called {\it acyclic}, an example being  Figure 3b.

{\bf 3B When all lines have quasi-isolated points.} We say that a PLS $(J,\Lambda)$ is a {\it QIMP}, if each $\ell\in\Lambda$ contains at least one {\it quasi-isolated} point, i.e. one which is on no other line. (The M in QIMP will be explained shortly.) For instance the PLS in Figure 3c (without $x$ and dashed line $\{4,x,10\}$) is a connected QIMP. Its line $\{1,6,5\}$ is interesting, while $\{5,7,8\}$ is not. In contrast,  the PLS in Figure 3b is {\it not} a QIMP since its 'horizontal' line has no quasi-isolated point. To complete the classification of lines in 3b: One line has one, three lines have two, and one line has three quasi-isolated points. In a QIMP every interesting line $\ell\in\Lambda$  has exactly  {\it one} quasi-isolated point. Furthermore, if $\ell$  occurs in any cycles $C$ and $C'$, it holds that:

(3.2)\quad ($C$-midpoint of $\ell)\ =$ ($C'$-midpoint of $\ell)\ =$ (unique quasi-isolated point of $\ell$)

{\bf 3B1 The structure of QIMPes.} All graphs appearing in this paper are assumed to be simple, i.e without multiple edges and loops. Here comes an easy 'structure theorem' for connected QIMPes. Start with any  connected graph $BP$, call it the {\it blueprint} graph, and 'plot' one new point on each edge. This obviously yields a connected QIMP $(J,\Lambda(BP))$, see Figures 4a, 4b. 
Conversely, let $(J,\Lambda)$ be {\it any} connected QIMP. Fix any $\ell=\{p,q,r\}$ in $\Lambda$ and let $q$ only be incident with $\ell$. If $(J_0,\Lambda_0)$ is defined by $J_0:=J\setminus\{q\}$ and $\Lambda_0:=\Lambda\setminus\{\ell\}$ then $(J_0,\Lambda_0)$ remains a QIMP.
Obviously $(J_0,\Lambda_0)$ has one or two connected components. We only follow up the latter possibility and leave the other to the reader. Namely, by induction the two components are induced by connected blueprint graphs having (say) $a$ and $b$ edges respectively. It follows that
$(J,\Lambda)$ is induced by a connected blueprint graph with $a+b+1$ edges.

We define the {\it midpoints of a QIMP} as the unique quasi-isolated points of its interesting lines.  Hence the acronym QIMP=quasi-isolated midpoints. The other points on the interesting lines are the {\it junctions of the QIMP}. For instance in the QIMP of Figure 4b each point is either a midpoint or a junction. That happens iff the QIMP is {\it 2-connected}, i.e. its blueprint graph is 2-connected in the usual sense. The QIMP in Figure 3c is not 2-connected. It has junctions 1,3,5, midpoints 2,4,6, but 7,8,9,10 are neither. See also 8A.

{\bf 3C When all interesting lines have unique midpoints.} We say the PLS $(J,\Lambda)$ is a {\it UMP} if for each interesting line $\ell\in\Lambda$  the first '=' in (3.2) takes place. Thus
UMPes generalize QIMPes in that again interesting lines have unique midpoints 'per se' but they need not be quasi-isolated.
 While acyclic $\not\Rightarrow$ QIMP, notice that acyclic $\Rightarrow$ UMP, vacuously by the lack of $C$-lines. 

For instance, $(J_1,\Lambda_1)$ in Figure 1c is not a UMP: If $C=(1',3',5',7')$ and $C'=(1',2',5',7')$ and $\ell=\{1',2',3'\}$ then the $C$-midpoint of $\ell$ is $2'$ but the $C'$-midpoint of $\ell$ is $3'$. Similarly the PLS in Figure 3a is not a UMP. 

\includegraphics[scale=0.52]{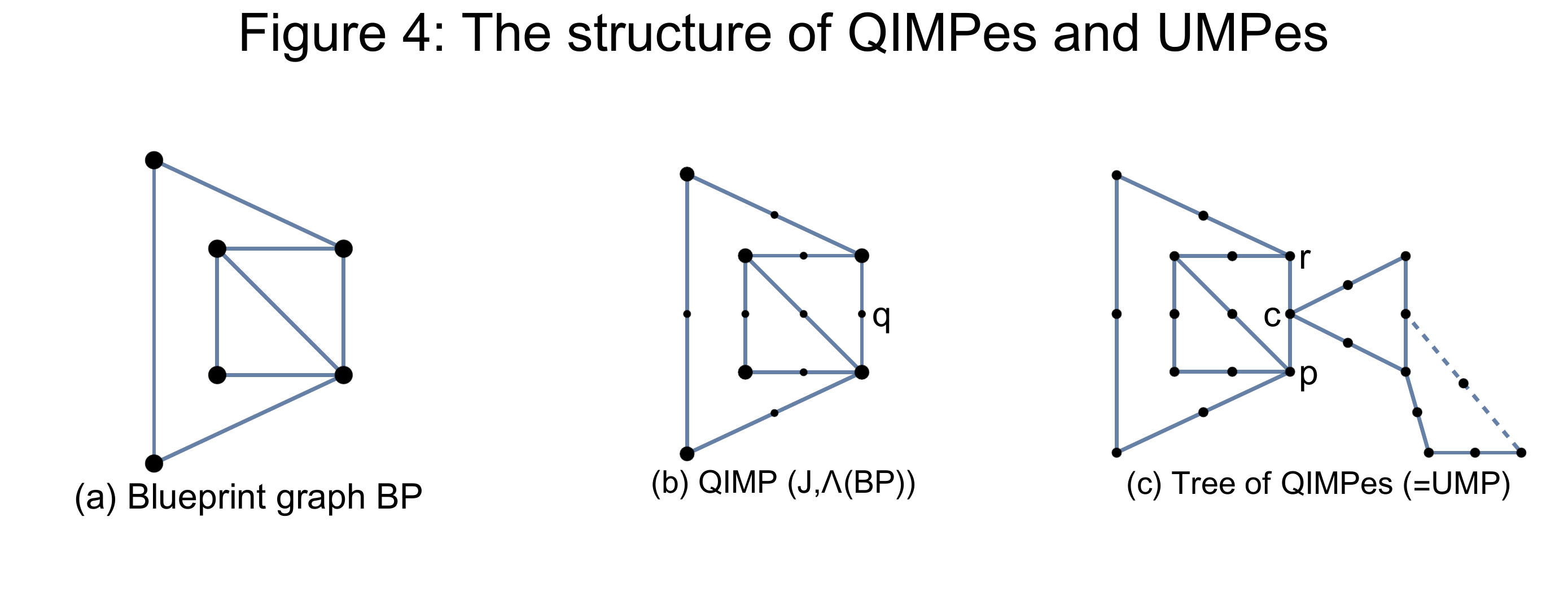}

 {\bf 3C1 The structure of UMPes.}  Let $(J^1,\Lambda^1)$ to $(J^t,\Lambda^t)$ be PLSes such that $J^2\cap J^1=\{p_2\},\ J^3\cap(J^1\cup J^2)=\{p_3\}$, and so on until $J^t\cap(J^1\cup\cdots\cup J^{t-1})=\{p_t\}$. Here the points $p_i$ need not be distinct. Putting $J:=\bigcup_{i=1}^t J^i$ and $\Lambda:=\bigcup_{i=1}^t \Lambda^i$ it is clear that $(J,\Lambda)$ is again a PLS. We call it a {\it tree} of the PLSes $(J^i,\Lambda^i)$. Using induction on $i$ one sees that each cycle $C$ of $(J,\Lambda)$ must be such that for some fixed index $j\in [12]$ all $C$-lines  are contained in $\Lambda^j$. It follows that each tree of QIMPes is a UMP. The converse holds as well [15,Lemma 13]:

{\bf Lemma 3.2:} The PLS $(J,\Lambda)$ is a UMP iff it is a tree of QIMPes.

For instance, gluing the midpoint $q$ of the QIMP in Figure 4b with the junction $1$ of the QIMP in Figure 3c yields the tree of QIMPes in Figure 4c. In this UMP the coalesced point $c$ is neither junction nor midpoint 'per se'. The tree structure in Lemma 3.2 is not at all unique. For instance the UMP in Figure 3b, even when shrunk to its largest connected component, is a tree of QIMPes in various ways. See also 8A.

{\bf 3D Benign midpoint-links.}  Let $(J,\Lambda)$ be a PLS containing a cycle $C=(p_1,\ldots,p_n)$, i.e. with ($C$-)junctions $p_i$ and ($C$-)midpoints $q_i\in [p_i,p_{i+1}]$. A path between two junctions $p_i$ and $p_j$ is harmless, whereas paths between $q_i,\ q_j$ or between $q_i,\ p_j$ will pose problems. Specifically,  a {\it type {\bf 1} midpoint-link} (of $C$) is a path between {\bf one} midpoint and a junction, i.e. $P=[q_i,\ldots,p_j]$ with $P^*\cap C^*=\{q_i,p_j\}$. And a {\it type {\bf 2} midpoint-link}  is a path $P$ between {\bf two} midpoints, i.e. $P=[q_i,\ldots,q_j]$ with $P^*\cap C^*=\{q_i,q_j\}$. We say $C$ has a {\it midpoint-link}, if the type doesn't matter.

\includegraphics[scale=0.4]{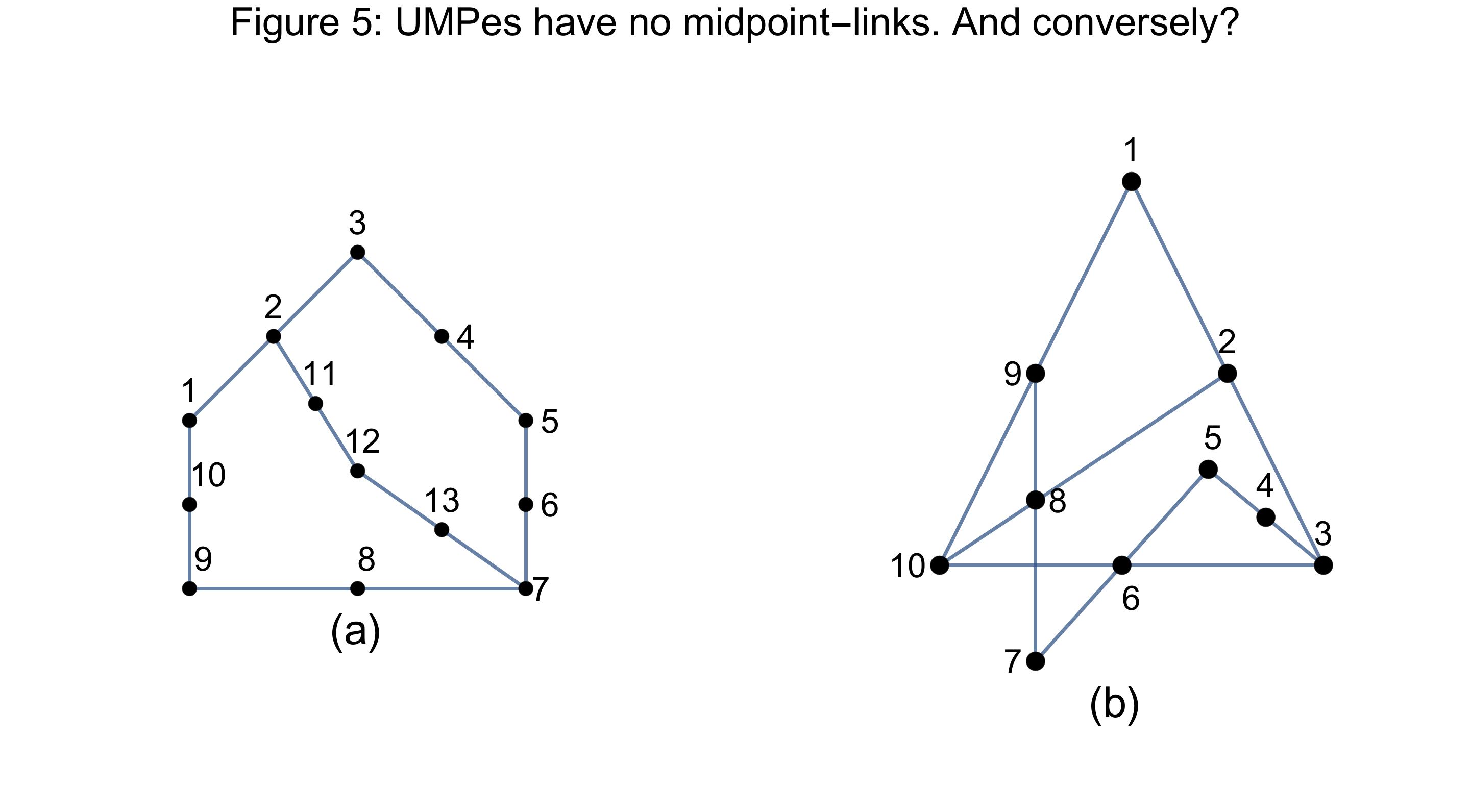}

The PLS $(J,\Lambda)$ in Figure 5a has a cycle $C_1=(1,3,5,7,9)$ with  a type 1 midpoint-link
$[2,12,7]$. Hence $C_1$ and $C_2=(2,3,5,7,12)$ have the common line $\ell=\{1,2,3\}$. Since the $C_1$-midpoint of $\ell$ is 2, whereas the $C_2$-midpoint is 1, our $(J,\Lambda)$ is not a UMP. A similar argument shows that also PLSes having  cycles with a type 2 midpoint-link are not UMPes.

By definition a PLS is a {\it NMPL} if its cycles have no midpoint-links. We just saw that UMP $\Rightarrow$ NMPL. Is it true that conversely 'not UMP $\Rightarrow$ not NMPL' ? For starters, 'not UMP' implies the existence of at least one {\it critical} cycle $C$ in the sense that $C$ has a line which, with respect to some other cycle, has another midpoint. By merely looking at the PLS $(J,\Lambda)$ in Figure 5a one could be misled to think that 'not UMP $\Rightarrow$ not NMPL' because in $(J,\Lambda)$ each critical cycle has {\it itself} a midpoint-link. (That {\it all} cycles in $(J,\Lambda)$ happen to be critical, is irrelevant.) Unfortunately matters are more intricate. One verifies that the PLS in Figure 5b has a critical cycle $(1,3,5,7,9)$ (e.g. in view of line $\{1,2,3\}$), yet $C$ itself has no midpoint-links. The fact that {\it other} cycles in 5b have midpoint-links only goes to show that 'not UMP $\Rightarrow$ not NMPL', if true at all, has a more subtle proof than 'UMP $\Rightarrow$ NMPL'.

{\bf Definition 3.3.} A type 1 midpoint-link $P$ of a cycle $C$ is {\it benign} if the $C$-midpoint and $C$-junction connected by $P$ are as close as they can get, i.e. they sit on the {\it same} line of $C$. Likewise, a type 2 midpoint-link of $C$ is {\it benign} if the two connected $C$-midpoints sit as close as they can get, i.e. they sit on {\it intersecting} lines of $C$. $\Box$

{\bf Example 3.4.} The cycle $(1,3,5)$ in Figure 6a has the benign type 1 midpoint-link $[6,8,10,5]$.
Similarly $C=(1,6,7,13)$ is a cycle whose $C$-midpoints of $[1,6]$ and $[6,7]$ are $5$ and $8$ respectively. Hence $[1,3,5]$ is a type 1 benign midpoint-link of $C$.
The cycle $(5,6,8,10)$ of the PLS $(J_2,\Lambda_2)$ in Figure 6a has the benign type 2 midpoint-link $[1,13,7]$. $\Box$

{\bf Definition 3.5.} The PLS $(J,\Lambda)$ is said to be a {\it BMPL} if all occuring midpoint-links are benign midpoint-links. $\Box$ 

Having no midpoint-links at all, each NMPL vacuously is a BMPL. In particular 'UMP $\Rightarrow$ BMPL'. For instance $(J_1,\Lambda_1)$ in Figure 1c is not a BMPL.
					
	\includegraphics[scale=0.67]{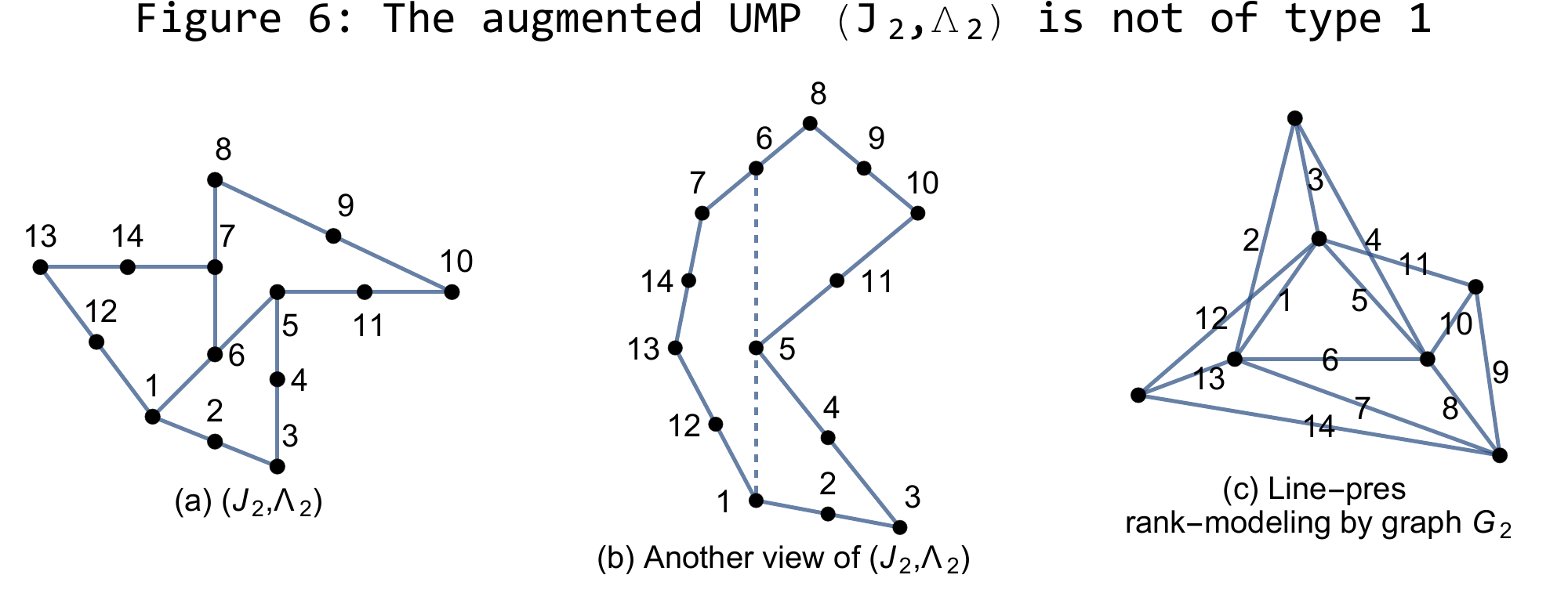}

					{\bf 3E Sparse PLSes.} We call a PLS $(J,\Lambda)$ {\it sparse} if there is a {\it testifying} ordering $(\ell_1,\ldots,\ell_t)$ of $\Lambda$ in the sense that $\ell_{i+1}\not\subseteq \ell_1\cup\ldots\cup\ell_i$ for all $1\le i<t$. Obviously 'QIMP $\Rightarrow$  sparse' since {\it any} ordering of $\Lambda$ will do. In view of Lemma 3.2, more generally 'UMP $\Rightarrow$  sparse'. Although $(J_1,\Lambda_1)$ in Figure 1c is not a UMP, it is sparse since

\quad $\ell_1=\{1',2',3'\},\ \ell_2=\{3',4',5'\},\ \ell_3=\{5',6',7'\},\ \ell_4=\{7',8',1'\},\ \ell_5=\{2',9',5'\}$

is a testifying ordering of $\Lambda_1$. 
In contrast, one checks that  $(J_4,\Lambda_4)$ in Figure 11a (Section 4) has 'too many lines' to be sparse.

{\bf 3F Augmented UMPes.}  We say that the PLS $(J',\Lambda')$ arises from the PLS $(J,\Lambda)$ by {\it adding a  path} if there is a path $P'=[x,z_1,...,z_s,y]\ (s\ge 0)$ of $(J',\Lambda')$ such that $J\cup P'^*=J'$ and $J\cap P'^*=\{x,y\}$ and 
					$\Lambda'=\Lambda\uplus\{[x,z_1],[z_1,z_2],...,[z_s,y]\}$. Clearly, when $(J,\Lambda)$ was sparse, then $(J',\Lambda')$ stays sparse. We say that $(J',\Lambda')$ arises from the PLS $(J,\Lambda)$ by {\it adding a benign midpoint-link} if  $P'$ is a benign  midpoint-link of some cycle  of $(J,\Lambda)$.
					
{\bf Example 3.6.} Letting $(J,\Lambda)$ be the cycle $C_1=(1,3,5)$ in Figure 6a, adding the benign (type 1) midpoint-link $[6,8,10,5]$ yields some PLS $(J',\Lambda')$. Then adding the benign (type 2) midpoint-link $[1,13,17]$ to the cycle $C_2=(5,6,8,10)$ of $(J',\Lambda')$ yields the PLS $(J_2,\Lambda_2)$. Apart from the obvious cycle $C_3=(1,13,7,6)$ there is a fourth cycle $C_4$ in $(J_2,\Lambda_2)$ which is shown in Figure 6b. One can show that $C_1$ to $C_4$ are the only cycles in $(J_2,\Lambda_2)$ and that their various midpoint-links are all benign. In particular notice that despite appearences the dashed line in Figure 6b is {\it not} a midpoint-link of $C_4$ since $\ell\subseteq C_4^*$. It follows that $(J_2,\Lambda_2)$
is a BMPL. (For the time being ignore Figure 6c.) $\Box$

 {\bf Example 3.7.} As we shall now see, unfortunately piling up benign midpoint-links {\it doesn't  always maintain} the BMPL property. Let us first verify that the PLS $(J,\Lambda)$ in Figure 7a (without the dashed line and its quasi-isolated point) is a BMPL. 

\includegraphics[scale=0.67]{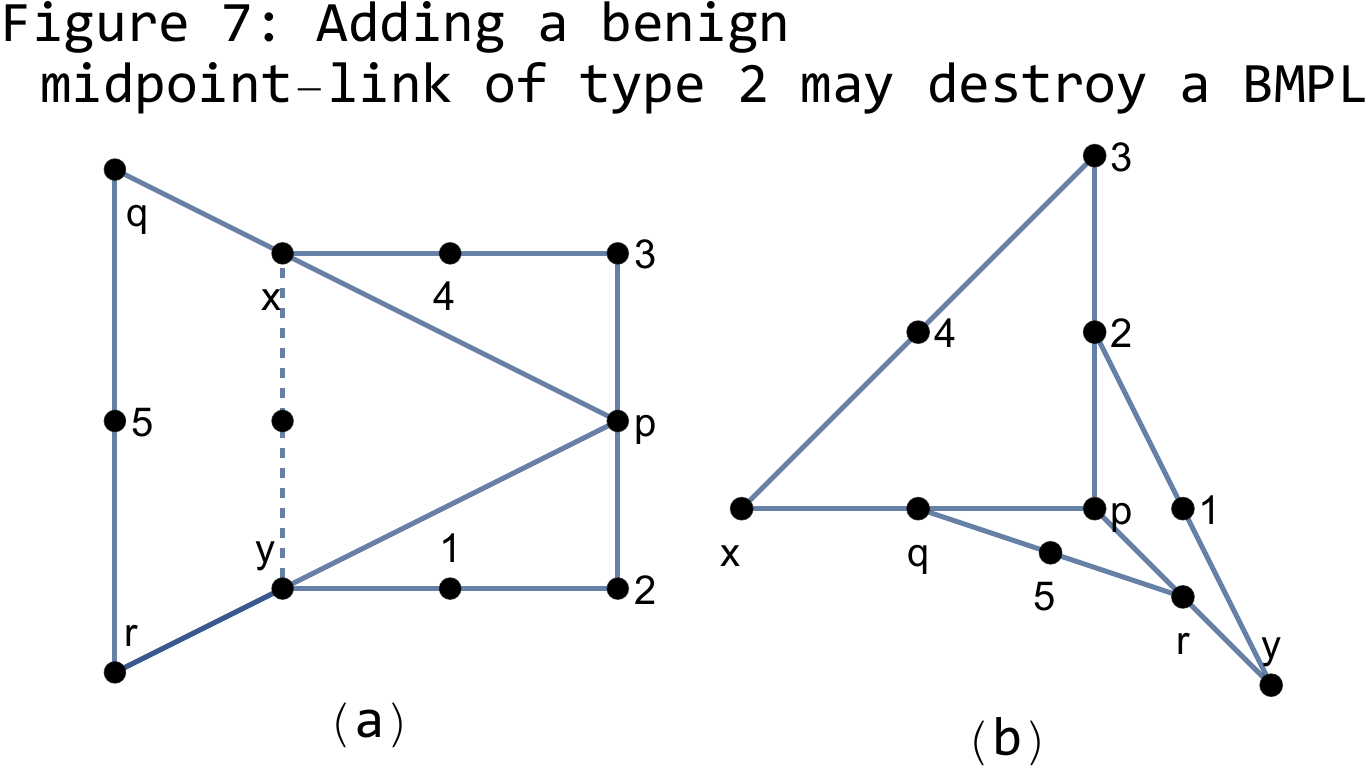}

Let $C$ be a cycle in $(J,\Lambda)$. Then at least one of the lines $[p,q],[q,r],[r,p]$ is a line of $C$. If at least two belong to $C$ then $C=(p,q,r)$, which has two benign midpoint-links $[x,3,p]$ and $[y,2,p]$. If exactly one of above-mentioned lines belongs to $C$ then this line is either $[p,q]$ or $[p,r]$. By symmetry we can assume it is $[p,q]$. This forces $C=(p,3,x)$, which has the benign midpoint-links $[p,y,2]$ and $[p,r,q]$, see Figure 7b. We thus verified that $(J,\Lambda)$ is a BMPL. But adding the (dashed) benign midpoint-link $[x,y]$ to the cycle $(p,q,r)$ produces a PLS $(J',\Lambda')$ which has a cycle $C'=(x,y,2,3)$ that features the non-benign midpoint-link $[x,p]$, see Figure 7a. $\Box$

 However, things improve when we stick to benign midpoint-links of {\bf type 1}: 

{\bf Lemma 3.8:}  Let the PLS $(J',\Lambda')$ arise from the BMPL $(J,\Lambda)$ by piling up benign type 1 midpoint-links.									
	\begin{enumerate}		
		\item[{(a)}]Then $(J',\Lambda')$ remains a BMPL.
			\item[{(b)}] If moreover $(J,\Lambda)$ has only midpoint-links of type 1, then so does $(J',\Lambda')$.
			\end{enumerate}
			
 If in contrast to Lemma 3.8 we allow piling up {\bf type 2} midpoint-links, this may (Fig 6a), or may not (Fig 7a), result in a BMPL. 

{\bf Definition 3.9.} An {\it augmented UMP} is a PLS that can be obtained (possibly in various ways) from a UMP by iteratively adding benign midpoint-links. The augmented UMP is {\it of type 1} if only benign midpoint-links of type 1 were added. $\Box$

Hence in a UMP  of type 1  {\it all} its midpoint-links are of type 1 (Lemma 3.8).		

{\bf Example 3.10.} Unfortunately not every BMPL is an augmented UMP. To witness, consider the PLS $(J_3,\Lambda_3)$ in Figure 9a. Notice that both $7$ and $8$ are quasi-isolated, i.e. belong to unique lines $\ell_1=\{2,8,9\}$ and $\ell_2=\{6,7,10\}$ respectively, whereas each other point $p\in J_3$ belongs to exactly two lines. Classifying the cycles $C$ of $(J_3,\Lambda_3)$ according to whether $C$ contains only $\ell_1$, or only $\ell_2$, or both of them, or none, yields the four cycles standing out in (a),(b),(c),(d) of Figure 9. All midpoint-links $P$ with respect to one of these cycles are benign. Specifically the 'square' $C$ in (a) has $P=[1,5]$, the 'square' in (b) has $P=[1,3]$, the 'pentagon' in (c) has no $P$, and the 'triangle' in (d) has $P=[4,10,6]$ and $P=[2,9,4]$. It follows that $(J_3,\Lambda_3)$ is a BMPL.	Yet these four benign midpoint-links are all {\it interfered with}; e.g. in (d) the midpoint $10$ of line $[9,4]$ in $P=[2,9,4]$ belongs to $\ell_2$. In other words, none of these $P's$ can be the last added path in an augmented UMP. $\Box$

\includegraphics[scale=0.45]{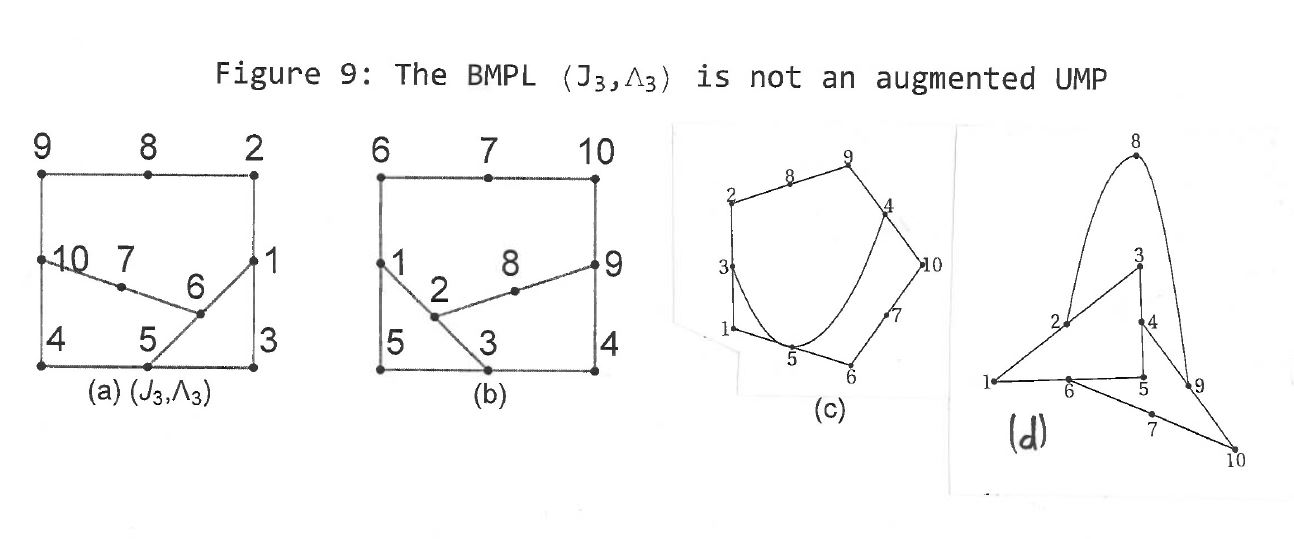}

{\bf Question 1:} As seen in Examples 3.7 and 3.10, neither is the class of augmented UMPes contained in the class of BMPLes, nor the other way round. What is the precise relation between these classes? And, harking back to 3D, is each NMPL a UMP?

{\bf 3G The PLS-rank.} 

{\bf Definition 3.11.} The {\it PLS-rank} (or simply {\it rank}) of $(J,\Lambda)$ is

(3.3) \qquad $rk(J,\Lambda):= |J|-|\Lambda|.$ $\Box$

 For instance $rk(J_1,\Lambda_1)=9-5=4$. In 8B we explain why this definition amounts to the seemingly quite different definition given in [15].
Let $(J_1',\Lambda_1')$ to $(J_c',\Lambda_c')$ be the connected components of the PLS $(J',\Lambda')$. Because of 

 $rk(J',\Lambda') = \ds\sum_{i=1}^c |J_i'|\ -\  \ds\sum_{i=1}^c |\Lambda_i'|\  =\  \ds\sum_{i=1}^c (|J_i'|-|\Lambda_i'|)\  =\ds\sum_{i=1}^c rk(J_i',\Lambda_i')$
					
		we will usually restrict ourselves to connected PLSes $(J,\Lambda)$ (with $\Lambda\neq\emptyset$). Then $\Lambda$  can be ordered in such a way $(\ell_1,\ldots,\ell_t)$ that $S_i:=(\ell_1\cup\ldots\cup\ell_i)\cap \ell_{i+1}\neq\emptyset$ for all $1\le i<t$. Upon adding $\ell_{i+1}$ the rank $rk$ of the PLS so far by $(3.3)$ changes in one of three ways:
		
	(3.4)\qquad If $|S_i|=1$ then $rk$ increases to $rk+(2-1)=rk+1$.
	
		(3.5)\qquad If $|S_i|=2$ then $rk$ remains $rk+(1-1)=rk$.
					
			(3.6)\qquad If $|S_i|=3$ then $rk$ {\it decreases} to $rk+(0-1)=rk-1$.

			Observe that the testifying orderings of $\Lambda$ in 3E are exactly the orderings that
					 avoid the rank-decreasing case (3.6).

	\section{Matroids and partial linear spaces}

	Recall from the Introduction that Section 4 links matroids [10] and PLSes, Section 5 links matroids and modular lattices, and Section 6 puts the pieces together. As to the finer structure of Section 4, see the summary given in 1B. We mention in addition that PLSes which can be modeled by a graph in a cycle-preserving manner are necessarily BMPLes (Lemma 4.11). A dual kind of property is circuit-friendliness.
	
	 {\bf 4A  Matroids that line-preserving model PLSes.} The definition below applies to arbitrary matroids but since we are mainly concerned with binary matroids,  recall that  $M(E)$ is {\it binary}  if the elements in $E$ can be matched with suitable vectors in some $GF(2)$-vector space $GF(2)^n$ such that dependency in $M(E)$ amounts to the linear dependency of the associated vectors. 
	
	{\bf Definition 4.1.} Let $M=M(E)$ be a simple [10] matroid with universe $E$. We say that a bijection  $\psi:J\to E$ {\it line-preserving (line-pres) models} the PLS $(J,\Lambda)$ if 
	
	(4.1)\qquad $\psi(\ell)$ is dependent in $M$ for all $\ell\in\Lambda$ {\it (dependence condition of the first kind)}. $\Box$

	{\bf Example 4.2.} The PLS in Figure 1c is line-pres modeled by a binary matroid as shown in Figure 10a. To unclutter notation the explicite mention of $\psi$ will often be omitted. If we assume that $a,b,c,d$ are independent (so $n\ge 4$) then the depicted 9-element matroid $M$ is  simple. 
	In contrast we claim that there is no binary matroid line-pres modeling the PLS $(J_3,\Lambda_3)$ from Figure 9a. By way of contradiction, assume the four 'corners' are $\psi$-labeled by some vectors $a,\ b,\ c,\ d$  in $GF(2)^n$, see Figure 10b. Because of (4.1) the midpoints of the cycle with junctions  $a,\ b,\ c,\ d$ must carry the labels $a+b,\ b+c,\ c+d$ and $a+d$. This forces a certain point being labeled by $(a+b)+(a+d)=b+d$, which in turn forces another being labeled by $(b+d)+(c+d)=b+c$. But now two points are labeled $b+c$, contradicting the simplicity of  $M(E)$. Similar reasoning shows (Figure 10c) that the PLS in Figure 7a cannot be line-pres modeled by a binary matroid. $\Box$

	\includegraphics[scale=0.72]{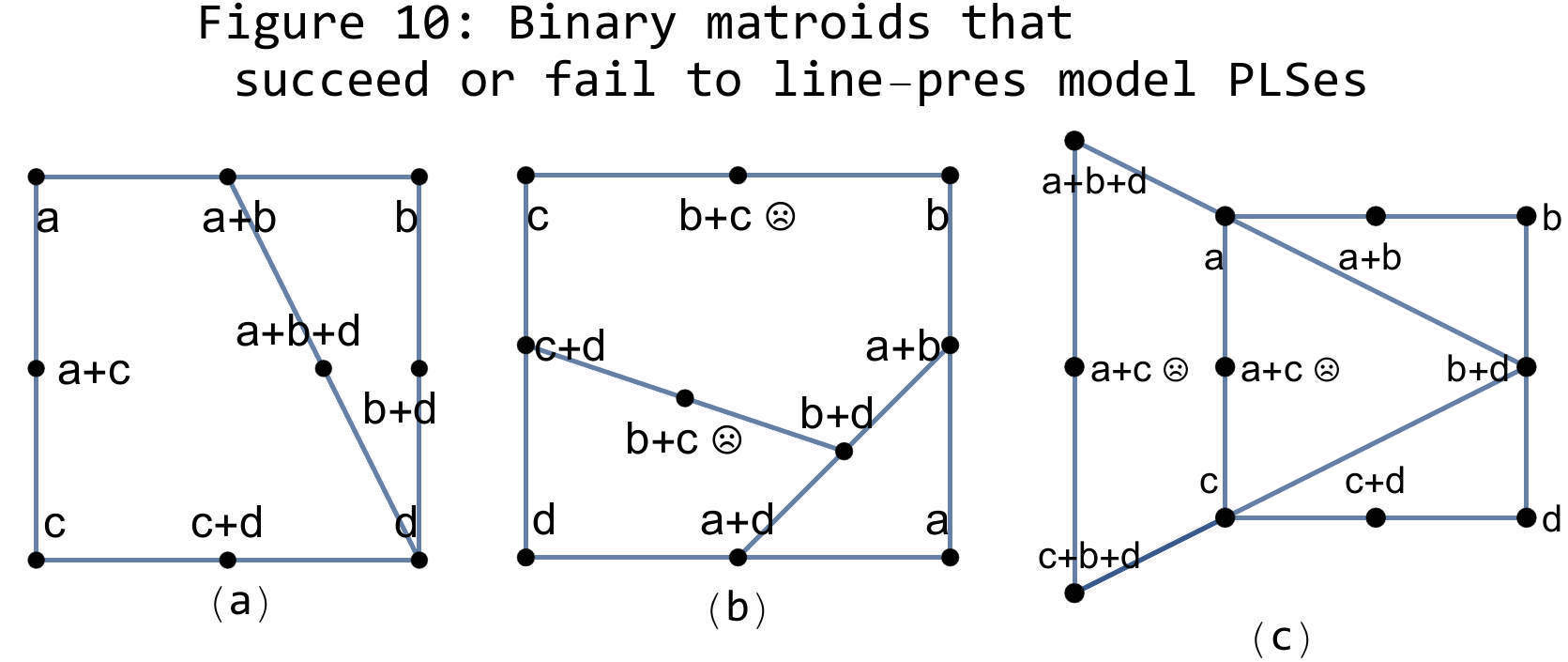}

	{\bf 4B Matroids that rank-model PLSes.}  Figure 11a shows a PLS $(J_4,\Lambda_4)$ with $J_4=\{a,b,c,p,q,r,x\}$ and $|\Lambda_4|=7$, as well as a binary matroid $M(E)$ with $E=GF(2)^3$ that line-pres models it. Thus say $x$ maps to $(0,1,1)$. (Figures 11b and 11c will only be relevant in Section 8F.)
	
	\includegraphics[scale=0.46]{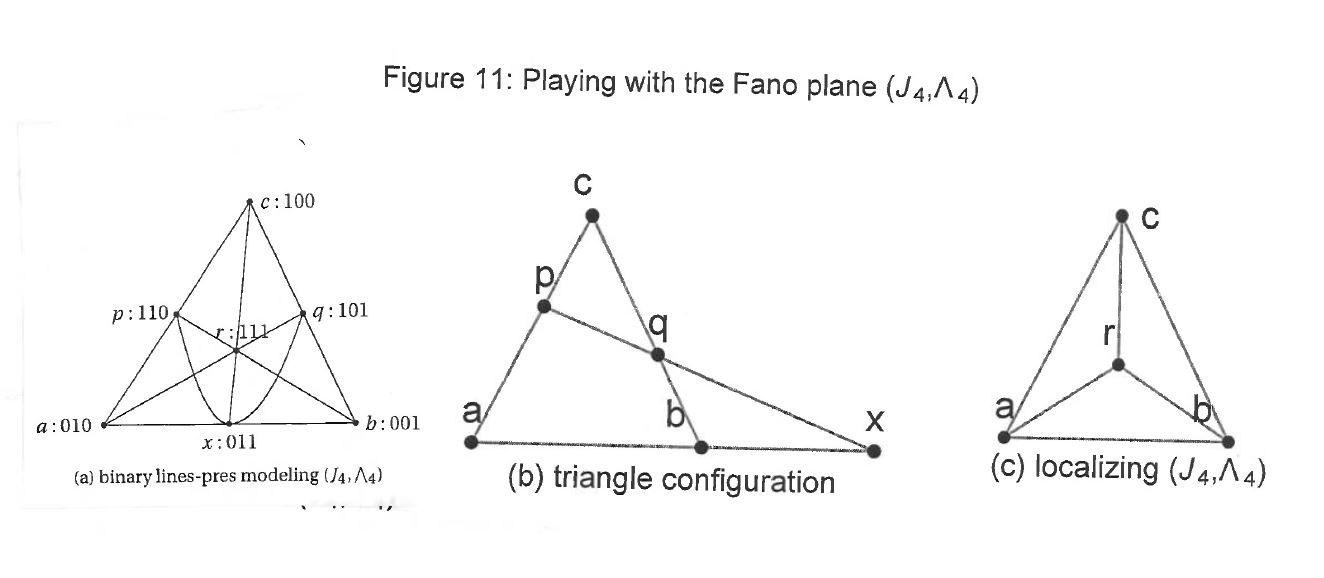}

	We have $rk(J_4,\Lambda_4)=7-7=0<3=mrk(E)$. 
	Actually the converse inequality will occur more often:
	
	{\bf Lemma 4.3:} If $(J,\Lambda)$ is sparse, and line-pres modeled by a matroid $M(E)$, then $mrk(E)\le rk(J,\Lambda).$
	
	{\it Proof.} Let $(\ell_1,\ldots,\ell_t)$ be a testifying ordering of $\Lambda$. For $t=1$ the inequality above becomes $2\le 2$. If generally $\Lambda'=\{\ell_1,\ldots,\ell_i\}$ and $E'\subseteq E$ matches $J':=\ell_1\cup\ldots\cup \ell_i$ via $\psi$ then $mrk(E')\le rk(J',\Lambda')$ by induction. Upon adding $\ell_{i+1}=\{a,b,c\}$ only cases (3.6) or (3.5) occur by sparsity. In case (3.6) the PLS-rank increases by 1. As to the matroid-rank, say $a\in E'$ and $b,c\not\in E'$. By submodularity $mrk(E'\cup\{b,c\})-mrk(E')\le mrk(\{a,b,c\})-mrk(\{a\})=1$, and so $mrk$ increases by at most 1. In case (3.5), the PLS-rank remains the same. As to the matroid side, say $a,b\in E'$ and $c\not\in E'$. Since $c$ is in the closure of $E'$ by (4.1), also the matroid-rank remains the same.
	 Thus the inequality gets perpetuated. $\Box$

	{\bf Definition 4.4.} We say a matroid $M(E)$ (or graph $G=(V,E)$) {\it rank-models} the PLS $(J,\Lambda)$ if 
	
	(4.2)\qquad $mrk(E)=rk(J,\Lambda)\ $  ({\it rank-condition of the first kind}.) $\Box$
	
	 Recall that $M(E)$ in Figure 10a  line-pres models $(J_1,\Lambda_1)$.
	 Since $rk(J_1,\Lambda_1)=4=mrk(E)$,  the matroid $M(E)$ also rank-models $(J_1,\Lambda_1)$. Observe that no matroid can rank-model
	$(J_4,\Lambda_4)$ in view of $rk(J_4,\Lambda_4)=0$.

	{\bf 4C Graphic matroids.} In the sequel we often focus on the most natural subclass of binary matroids $M(E)$, i.e. the class of {\it graphic} matroids. Thus by definition $E$ is the edge set of a graph $G=(V,E)$ with vertex set $V$, and a subset of $E$ is dependent iff it contains the edge set of a circuit [10].  There is no harm identifying a circuit with its underlying edge set.
	Consequently if the graphic matroid $M(E)$ line-pres models $(J,\Lambda)$  then each line $\ell$ maps to a triangle of $G$. Instead of saying '$M(E)$ line-pres models $(J,\Lambda)$' we usually say the {\it graph $G$ line-pres models} $(J,\Lambda)$.
	
	{\bf 4C1 The naive algorithm.} Here comes a simple (exponential time) algortithm that decides whether or not any given PLS $(J,\Lambda)$ can line-pres be modeled by a suitable graph (that gets constructed in the process). It works as follows. Match each line $\ell\in\Lambda$ with an isolated triangle whose edges are labeled with the points of $\ell$. Start with any $p\in J$ and glue together all edges of triangles labeled by $p$. Call the obtained partial graph $G(p)$. Then pick any label $q\neq p$ occuring in $G(p)$ and repeat. This yields $G(p,q)$, and so on. The procedure either gets stuck or leads to the desired graph $G$.  

\includegraphics[scale=0.52]{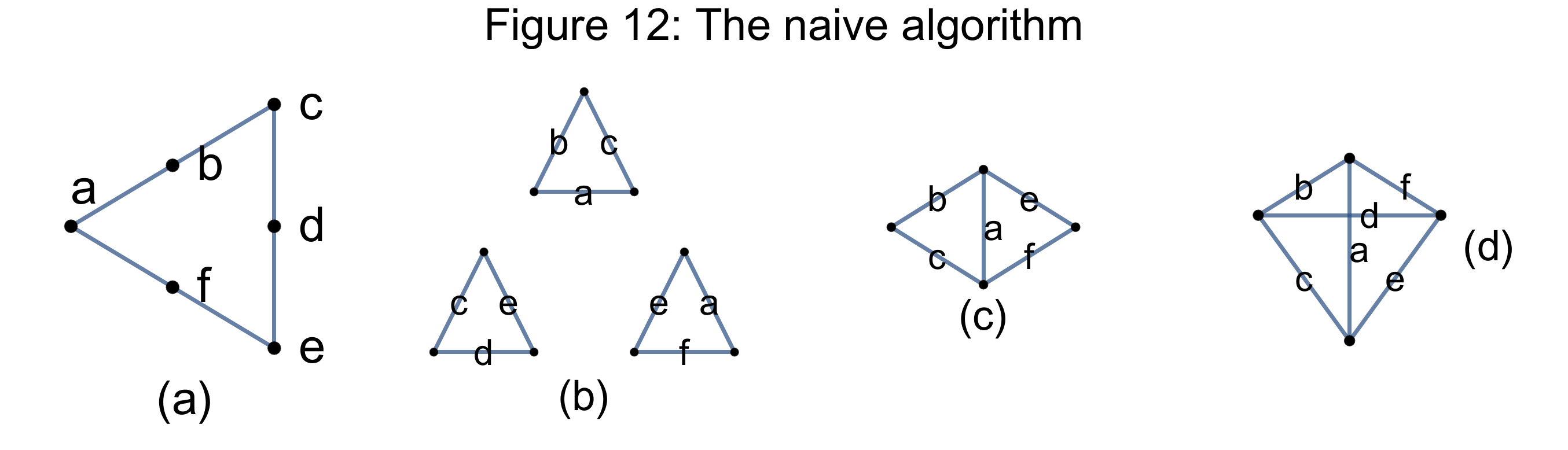}

For instance, the PLS in Figure 12a triggers the labeled triangles in Figure 12b. There are two triangles with an edge labeled $a$, and they can e.g. be glued to a graph $G(a)$ as in Figure 12c. But $G(a)$ is a dud since edges $c$ and $e$ are not incident. However, if one glues the two triangles as in Figure 12d (i.e. the triangle $\{a,e,f\}$ got {\it flipped}), then the edges $c$ and $e$ can be completed to a triangle $\{c,e,d\}$. The resulting graph $G$ line-pres models our PLS.

 In general, whenever $t$ triangles need to be glued along an edge, we do so in $2^t$ ways (flip, or not, each triangle), and only persue the partial graphs which are no duds. More clever ways to decide line-preseving modeling follow in due course.

	{\bf Question 2:}  What are necessary or sufficient conditions beyond the algorithm of 4C1 for a PLS to possess a line-preserving modeling graph? What about line-pres rank-modeling graphs?

	{\bf 4C2 The standard graph of a QIMP.} We now turn to QIMPes and line-pres rank-model them with graphs. Let $(J^1,\Lambda^1)$ be a connected QIMP.
	We may assume (see 8A for justification) that $(J^1,\Lambda^1)$ is 2-connected, i.e. each point is either a junction or a midpoint. Thus the set $\{p_1,\ldots, p_n\}$  of junctions can be identified with the vertex set of the blueprint graph, see 3B1. As in [15, p.218]
	we put $V^1=\{0,1,\ldots,n\}$ and construct a graph  $G^1=(V^1,E^1)$ that will line-pres model $(J^1,\Lambda^1)$. Before we define $\psi:J^1\ra E^1$ we hasten to stress that $G^1$ is not to be confused with the blueprint graph, e.g. $G^1$ has one vertex more. As visualized in the first {\it column} of Figure 13,
	each junction $p_i$ is mapped (by $\psi$) to the edge $\{0,i\}$.  Furthermore, if the line $[p_i,p_j]$ exists then its midpoint $q$ is mapped to $\{i,j\}$. In this way $(4.1)$ is satisfied.  (At this stage the dashed edge is not special.) Actually rank-modeling takes place as well: $rk(J^1,\Lambda^1)=6-3=|V^1|-1=mrk(E^1)$.

	The construction of $G^1$ generalizes to arbitrary QIMPes. The obtained line-pres rank-modeling graph, call it again $G^1$, is the {\it standard graph of the QIMP}. Its standout vertex is $0$: The $C$-junctions of any fixed cycle $C$ of the QIMP map to certain edges (call them 'spokes')  among the $n$ edges incident with $0$. Furthermore, the $C$-midpoints map to the 'rimes' of a {\it wheel} 
	$W$ having these spokes. Thus these transformations occur:
	
	$\begin{array}{ccc}
	\hbox{\it junctions of C} & \rightarrow & \hbox{\it spokes of W} \\
	\hbox{\it midpoints of C} & \rightarrow & \hbox{\it rims of W}
	\end{array}$
	
While $G^1$ may feature many wheels, all of them are {\it centered at $0$}. We shall call a wheel {\it degenerate} if it has $n=2$, spokes, i.e. if it is a triangle. A triangle is the only kind of wheel whose set of spokes is not uniquely determined.

Let us illustrate all of this once more. The second column in Figure 13 shows that the standard graph of the QIMP $(J^2,\Lambda^2)$ is $G^2=(V^2,E^2)$, which yields a rank-modeling since
	$rk(J^2,\Lambda^2)=5-2=|V^2|-1=mrk(E^2)$.

	{\bf 4C3 The standard graph of a UMP.} Generally, given a tree-PLS in the sense of 3C1, and modeling graphs of the component PLSes, these graphs can be combined to a modeling graph of the tree-PLS [15, Lemma 12]. We merely illustrate this  for the two component QIMPes in Figure 13 and their cooresponding standard graphs.  Merging $q$ of $(J^1,\Lambda^1)$ with $p_1'$ of $(J^2,\Lambda^2)$ yields the tree of QIMPes (=UMP) top right in Figure 13. Continuing our indexing scheme let's denote this UMP by $(J_5,\Lambda_5)$. For later use note that

	\begin{enumerate}
	\item [(4.3)]  $rk(J_5,\Lambda_5)=|J_5|-|\Lambda_5|=|J^1|+|J^2|-1-|\Lambda^1|-|\Lambda^2|=\, rk(J^1,\Lambda^1)+rk(J^2,\Lambda^2)-1.$
	%\item[{}] $=\, rk(J^1,\Lambda^1)+rk(J^2,\Lambda^2)-1.$
	\end{enumerate}

	\includegraphics[scale=0.35]{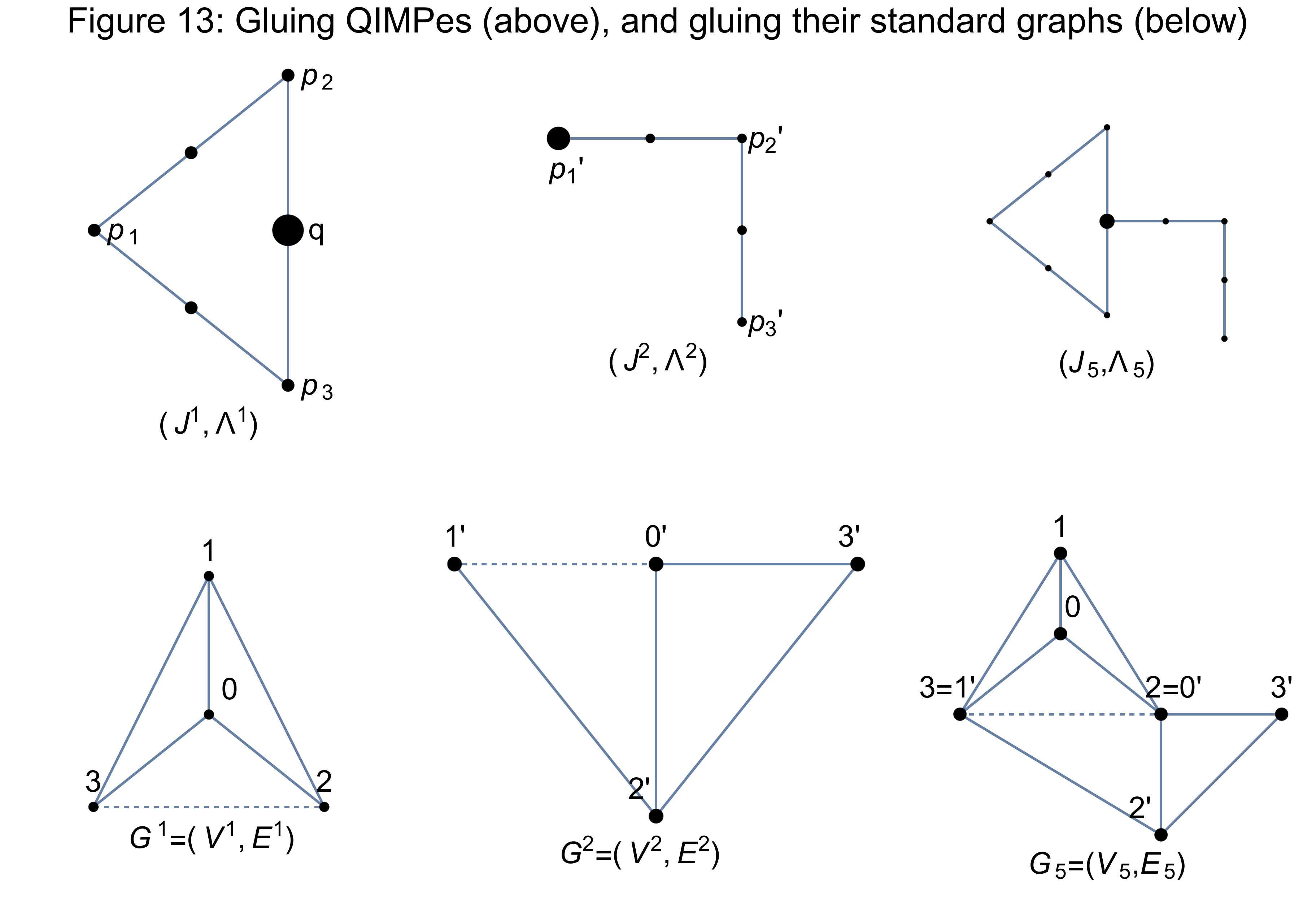}

	 If by mimicking the merging of the points $q$ and $p_1'$, we merge the corresponding (dashed) edges $\psi(q)=\{3,2\}$ of $G^1$ and $\psi(p_1')=\{1',0'\}$ of $G^2$, then we get a graph $G_5=(V_5,E_5)$. It evidently line-pres models $(J_5,\Lambda_5)$. In fact rank-modeling (4.2) carries over as well:

	\begin{enumerate}
	\item [(4.4)]
	$mrk(E_5)=|V_5|-1=|V^1|+|V^2|-2-1=mrk(E^1)+mrk(E^2)-1$
	\item[{}] $\,{(4.2)\atop =}\,rk(J^1,\Lambda^1)+rk(J^2,\Lambda^2)-1\,{(4.3)\atop =}\,rk(J_5,\Lambda_5)$
	\end{enumerate}
	
	In view of Lemma 3.2 we have thus sketched a proof of [15, Lemma 14(a)], which we restate as follows.
	
	{\bf Lemma 4.5:} Each UMP can be line-pres rank-modeled by a graph.
	
	We call the graph $G$ in Lemma 4.5 the {\it standard graph of the UMP}. Having grasped how  $G$ arises, it will henceforth be more convenient to {\it label the edges} rather than the vertices of $G$. In this way we can use the {\it same labels} as for the points of the UMP. Thus the third column of Figure 13 gives way to Figure 14.

	\includegraphics[scale=0.36]{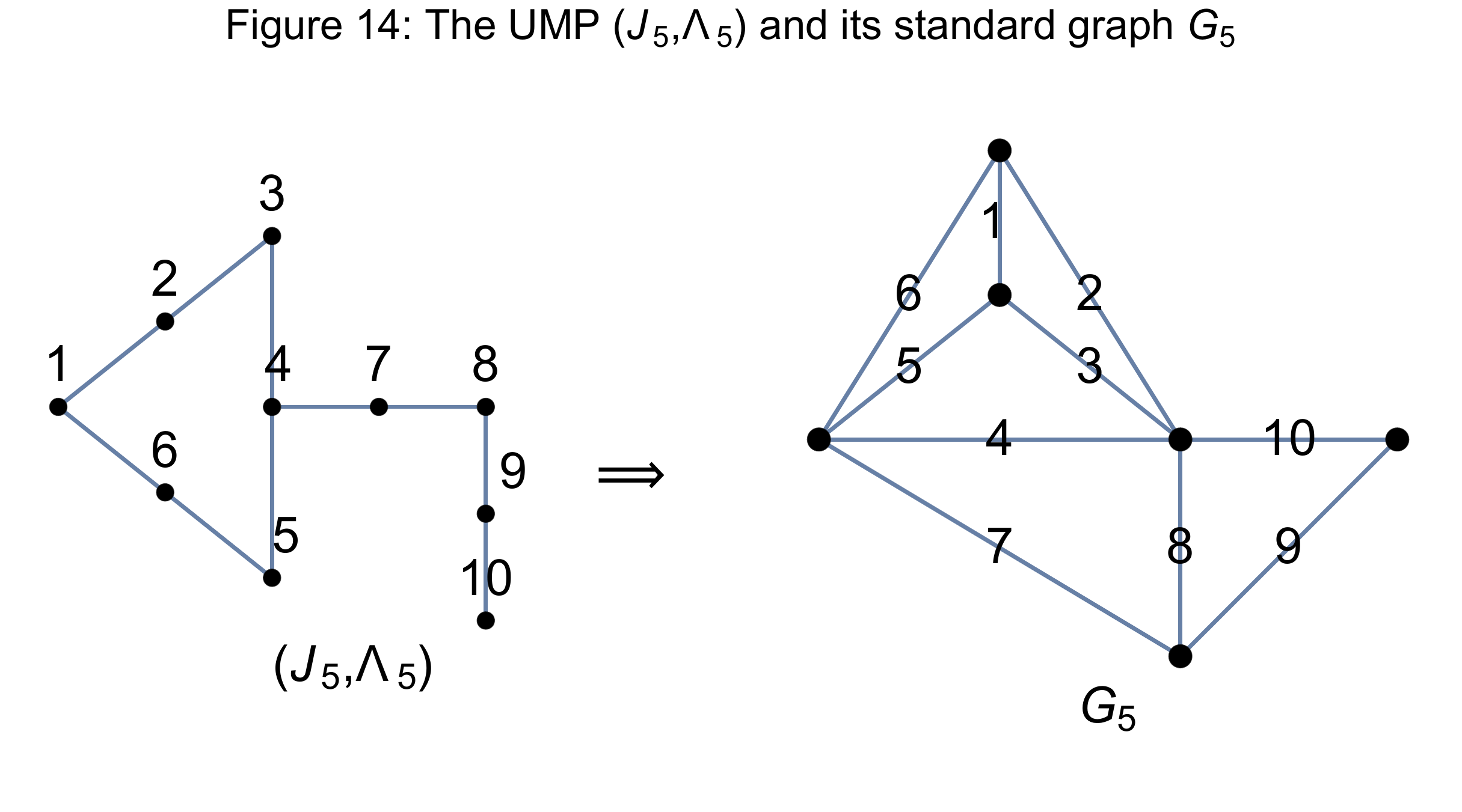}

	{\bf 4C4 Nonstandard modeling of cycle-PLSes.} Each cycle $C=(p_1,\ldots,p_n)$ in a PLS induces a {\it cycle-PLS} $(C^*,\Lambda^*)$ where $\Lambda^*:=\{[p_i,p_{i+1}]:\ 1\le i\le n\}$ and $C^*$ is (recall) the underlying set of $C$-junctions $p_i$ and $C$-midpoints $q_i\in[p_i,p_{i+1}]$. We say that $(C^*,\Lambda^*)$ {\it comes from} the cycle $C$. Since the number $n$ of $C$-junctions  equals $|\Lambda^*|$, and whence the number of $C$-midpoints, we have $|C^*|=n+n$, and so 
	
	(4.5)\quad $ rk(C^*,\Lambda^*)=|C^*|-|\Lambda^*|=|\Lambda^*|.$
	
	Each cycle-PLS is
	 a QIMP whose standard graph boils down to a single wheel $W=(V,E)$, see Figure 15 where $n=6$. Furthermore, it is easy to see that {\it each} line-pres modeling of a cycle-PLS  by a wheel must be the standard modeling, and so is rank-preserving.

	\includegraphics[scale=0.4]{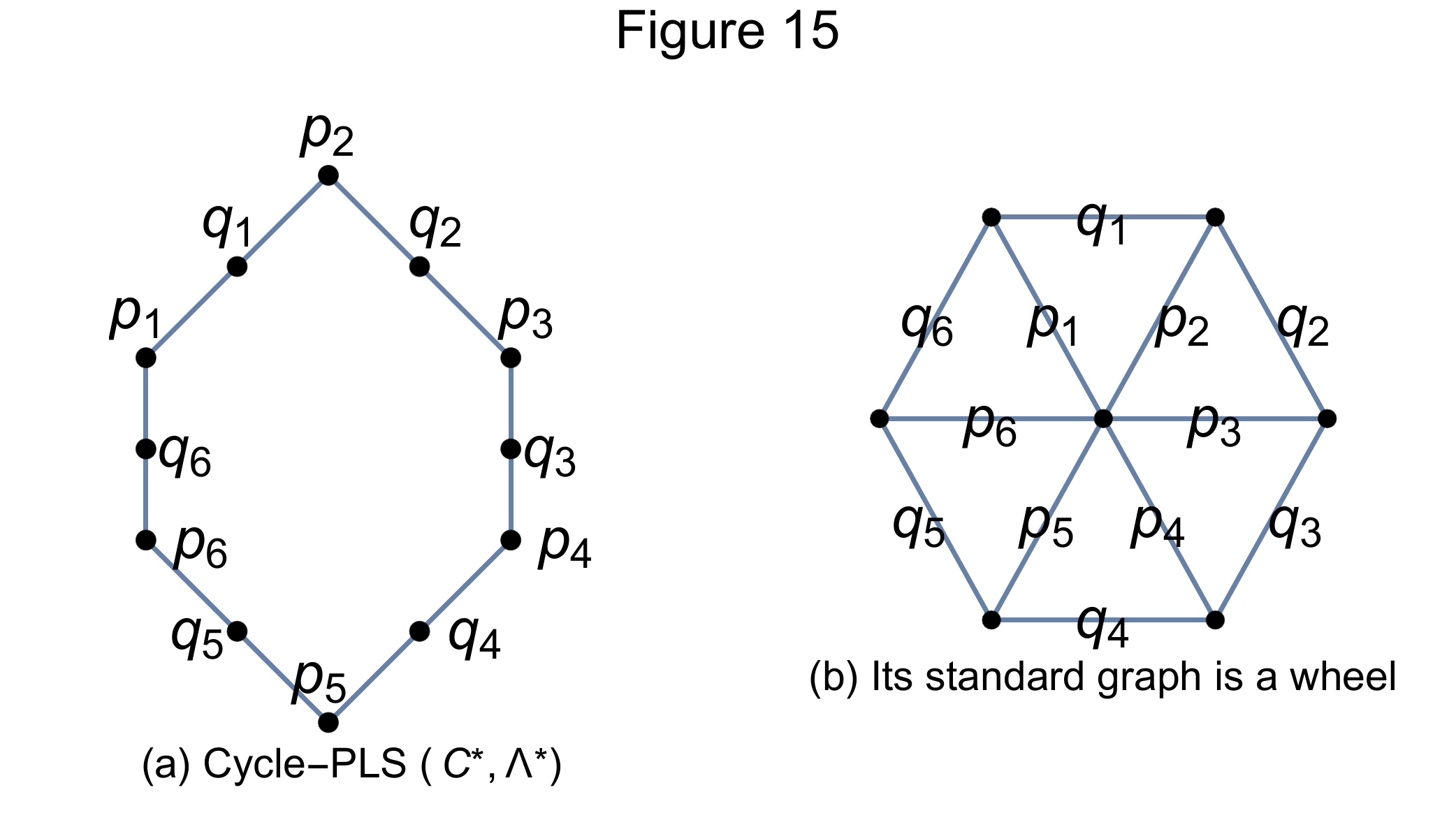}
	
	Let us investigate whether cycle-PLSes  can be line-pres modeled by {\it other} graphs $G=(V,E)$, and if so, must be rank-preserving. To begin with, since cycle-PLSes are sparse it follows from Lemma 4.3 that
	
	{\bf Corollary 4.6:} For each line-pres graph modeling $\psi:C^*\ra E$ one has $mrk(\psi(C^*))\le rk(C^*,\Lambda^*)$.

	Of course if $\psi$ rank-models $(C^*,\Lambda^*)$ then by definition $=$ takes place in Corollary 4.6.
	That $<$ can take place for merely line-pres modelings is witnessed by the cycle-PLS $(C^*,\Lambda^*)$ in Figure 16, which has a line-pres modeling by a {\it non-wheel} graph 
	$G=(V,E)$.  
	Indeed in view of (4.5) one has  $mrk(E)=|V|-1=4<5=|\Lambda^*|=rk(C^*,\Lambda^*)$.
	
	\includegraphics[scale=0.39]{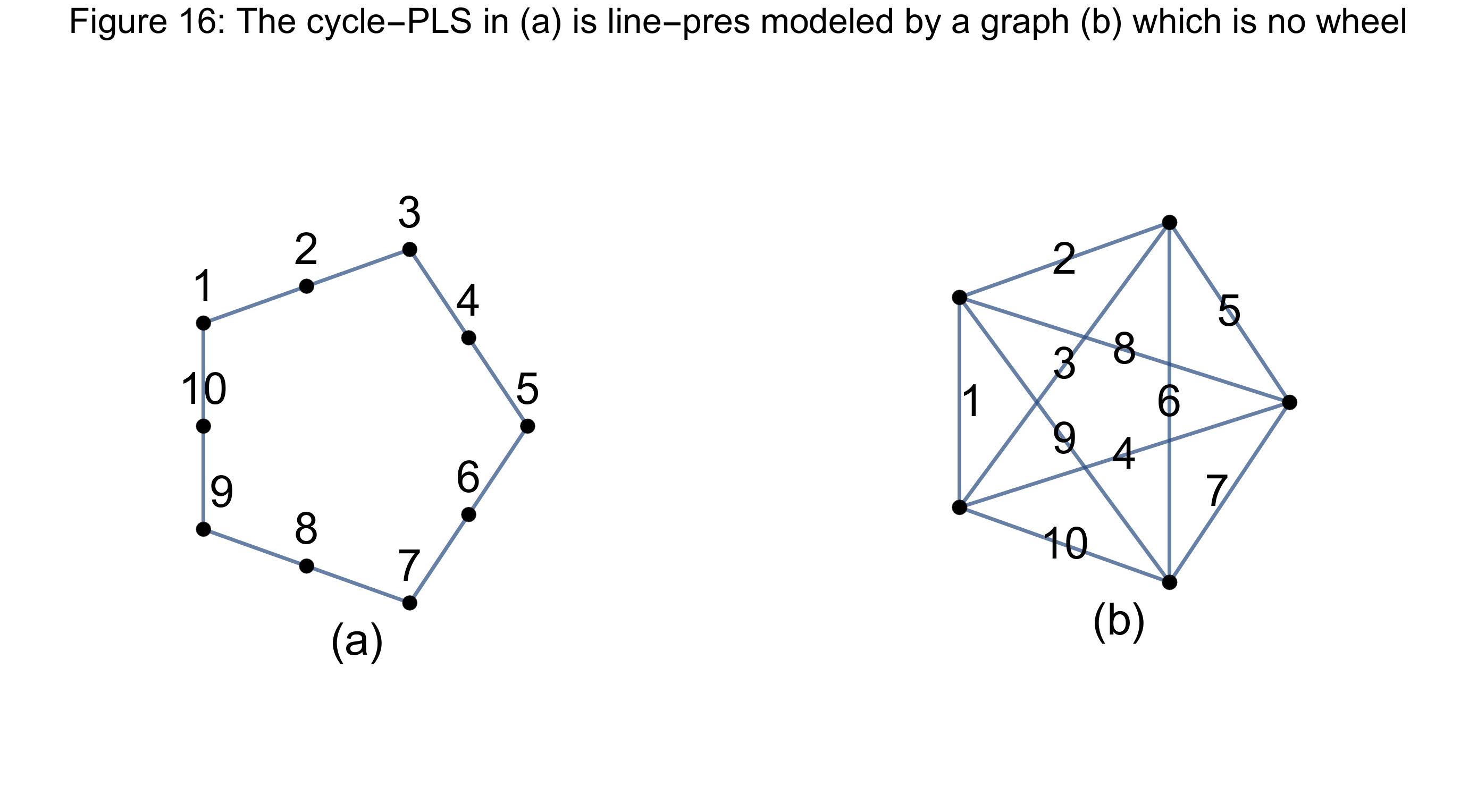}
	
	 The good news is, not only is each line-pres modeling of a cycle-PLS by a wheel automatically rank-preserving, the converse holds as well:

	{\bf Lemma 4.7:} Let the PLS $(J,\Lambda)$ be line-pres modeled by $G=(V,E)$ via $\psi:J\ra E$. If $(C^*,\Lambda^*)$ comes from a cycle $C$  of $(J,\Lambda)$ 
	and  $rk(C^*,\Lambda^*)=mrk(\psi(C^*))$, then $\psi(C^*)$ is a wheel.

The bad news is: When $G=(V,E)$ line-pres rank-models $(J,\Lambda)$, thus $rk(J,\Lambda)=mrk(E)$ by definition, this does not imply that  $rk(C^*,\Lambda^*)=mrk(\psi(C^*))$ for all cycles $C$ of $(J,\Lambda)$!
 
{\bf Example 4.8.} For instance $G_2=(V_2,E_2)$ in Figure 6c line-pres rank-models $(J_2,\Lambda_2)$ in Figure 6b since (4.1) is verified ad hoc (i.e. lines map to triangles), and $rk(J_2,\Lambda_2)=14-8=7-1=mrk(E_2)$. The large cycle $C$ in Figure 6b has $C^*=J_2$ (and so $\psi(C^*)=E_2$) but 
$(C^*,\Lambda^*)\not =(J_2,\Lambda_2)$ since $\Lambda_2\setminus\Lambda^*=\{\{1,5,6\}\}$. Hence 	
$rk(C^*,\Lambda^*)=|\Lambda^*|=7>6=mrk(\psi(C^*))$.	$\Box$
		
	{\bf 4C5 Cycle-preserving modeling of PLSes.} Example 4.8 motivates the following concept. A line-pres  graph $G$ of a PLS $(J,\Lambda)$ is {\it cycle-preserving (cycle-pres)} if the midpoints of each cycle $C$ of $(J,\Lambda)$ map to a circuit $\Gamma$ of $G$. This forces $\Gamma$ to be of a very specific shape. Namely, since each $\ell\in\Lambda$ maps to a triangle in $G$, a quick sketch confirms that the edges in $\Gamma$ are the rimes of a wheel whose spokes bijectively correspond to the junctions of $C$.
	Hence by Lemma 4.7 we conclude: 
	
	{\bf Corollary 4.9:} Let $\psi$ be a line-pres modeling of a PLS. Then $\psi$ is cycle-pres iff for all cycles one has $rk(C^*,\Lambda^*)=mrk(\psi(C^*))$.

	The glaring example of a cycle-preserving graph is, by its very definition, the standard graph of a QIMP.
	Since by Lemma 3.2 each cycle in a UMP is contained in one of its QIMP components we conclude:

			{\bf Corollary 4.10:} The standard graph of a UMP is cycle-preserving.
			
Recall from Subsection 3D that UMPes generalize to BMPLes. Here comes a kind of converse.

					{\bf Lemma 4.11:} If $(J,\Lambda)$ admits a cycle-pres modeling graph $G$ then $(J,\Lambda)$ is a BMPL.

					{\it Proof.} Consider a cycle $C$ of $(J,\Lambda)$, without much loss of generality let's take the one in Figure 15a. By cycle-preservation it is mapped onto the rimes of a wheel of $G$ (see 15b). By way of contradiction assume  $C$ had non-benign midpoint-links.

	{\it Case 1:} Suppose there is a non-benign type 1 midpoint-link, such as $(q_2,p_4)$ in Figure 18a.
					As shown in Figure 18b this yields a cycle $C'$ in $(J,\Lambda)$ with junctions $p_4,\ q_2$ (and $p_3$). By cycle-preservation $C'$ maps to a wheel of $G$, in such a way
					that $p_4,\ q_2$  are mapped to incident edges (being spokes).  But this contradicts Figure 15b where these edges are not incident.
					(Generally in each wheel $p_i$ is  incident with  $q_j$ only when $j\in\{i-1,i\}$.)

					{\it Case 2:} Suppose there is a non-benign type 2 midpoint-link, such as $(q_2,q,q_5)$ in Figure 18c. Then, as
					 illustrated in Fig. 18d, there is a cycle $C''$ in $(J,\Lambda)$ with junctions (among others) $q_2$ and $q_5$. As in case 1 they are mapped to incident edges of $G$. This contradicts Figure 15b where $q_2,\ q_5$ do not touch. 
					(See 8C for further remarks.)  $\Box$

					\includegraphics[scale=0.48]{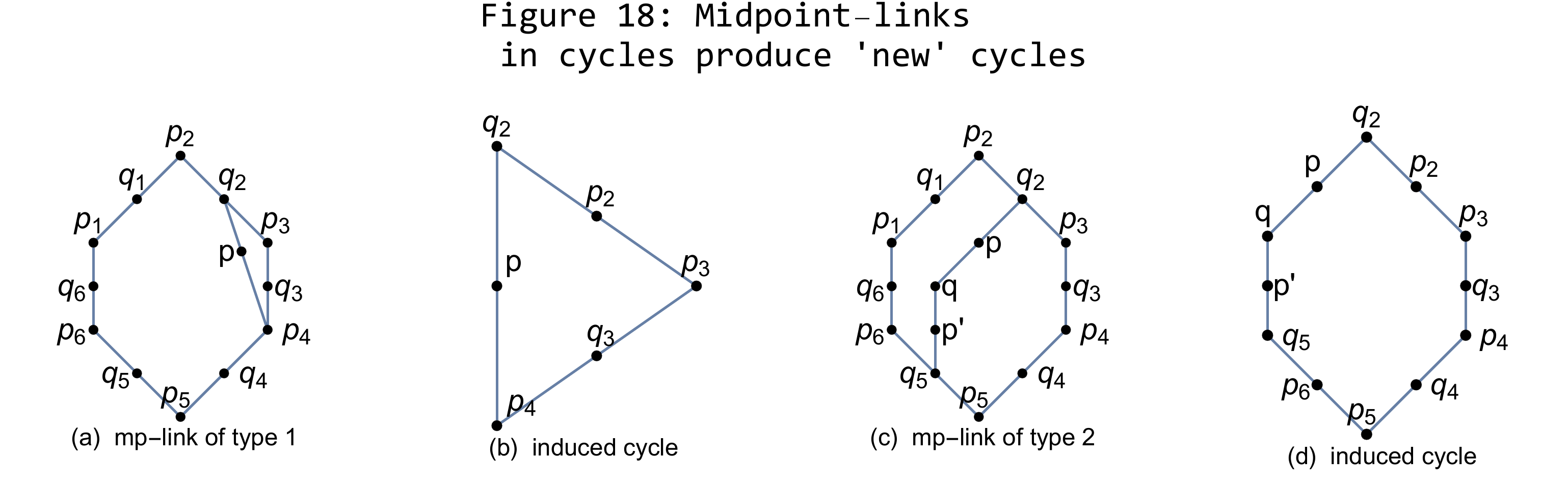}

				{\bf 4C6 Circuit-friendly modeling of PLSes.} Here comes a property which in a sense is the dual of cycle-preservation, and which in fact  will be more relevant than cycle-preservation in Section 6. We say that a graph $G=(V,E)$ models the PLS $(J,\Lambda)$ {\it circuit-friendly} if there is a bijection $\psi: J\ra E$ such that
				 each chordless circuit $\Gamma$ of $G$ either corresponds to a line (then necessarily $|\Gamma|=3$), or 
				to the midpoints of a cycle of $(J,\Lambda)$. Hence triangles in $G$ are chordless circuits with two options to behave under circuit-friendliness. For instance in Figure 14 the triangle $\{2,4,6\}$ of $G_5$  corresponds to the midpoints of a cycle of the UMP $(J_5,\Lambda_5)$, whereas the triangle $\{1,2,3\}$ matches a line. The circuit $\{4,7,9,10\}$ of $G_5$ is allowed to behave neither way since it is has the chord 8.

				Recall  that the standard graph $G$ of a QIMP is cycle-preserving 'by construction'. It is more subtle [15, Lemma 14(b)] that $G$ is circuit-friendly as well. This being established it follows (as in Corollary 4.10) from Lemma 3.2 that QIMPes generalize to UMPes:

				{\bf Lemma 4.12:} The standard graph of a UMP is circuit-friendly.

	To summarize 4C5 and 4C6, cycle-preservation pressurizes the {\it  cycles} $C\subseteq J$ in a PLS (with respect to a bijection $J\ra E$), while circuit-friendliness pressurizes the {\it  chordless circuits} $\Gamma\subseteq E$.
				
				{\bf 4C7 Modeling BMPLes with graphs.} Lemma 4.11 showed that under natural circumstances BMPLes are necessary. Are they sufficient? Unfortunately, the BMPL puzzles raised in Question 1 accumulate when it comes to being modeled by  graphs. To witness, {\it not every} BMPL $(J,\Lambda)$ can occur in Lemma 4.11 since $(J,\Lambda)$ may not have a line-pres modeling graph at all. Worse, the BMPL $(J_3,\Lambda_3)$ of Figure 9, that underlies Figure 10b, doesn't even admit a  binary line-pres modeling matroid. However this BMPL is
				not an augmented UMP. Do at least the latter possess line-pres rank-modeling graphs? In order to give  partial answers in Lemma 4.14 (b),(c) we need a few preliminaries.
				
Recall from 3F the concept of adding a path $[p_1,\ldots,p_n]$ to a PLS. The $n=2$ scenario boils down to having two points $p_1,\ p_2$ of a PLS $(J,\Lambda)$ that is line-pres rank-modeled by a graph $G=(V,E)$ (so $|J|-|\Lambda|=|V|-1$), and to asking whether $p_1,\ p_2$ can be connected by a line $\ell$ while extending $G$ accordingly. Specifically, let $q\not\in J$ be a new point and extend $(J,\Lambda)$ to $(J',\Lambda')$ where $J'=J\cup\{q\}$ and $\Lambda'=\Lambda\cup\{\ell\}$ with $\ell:=\{p_1,p_2,q\}$. If the edges $\ol{p_1},\ \ol{p_2}\in E$ coupled to $p_1,\ p_2$ happen to be {\it incident}, then we can extend $G$ to a graph $G'=(V',E')$ with no extra vertex but one extra edge $\ol{q}$ that yields a triangle $\{\ol{p_1},\ol{p_2},\ol{q}\}$. Thus $G'$ remains line-preserving, and rank-preservation is perpetuated as well:

$rk(J',\Lambda')=|J'|-|\Lambda'|=(|J|+1)-(|\Lambda|+1)=|V|-1=|V'|-1=mrk(E').$

All of this can be extended to $n>2$ if one postulates again that $\ol{p_1},\ \ol{p_n}$ be incident (Lemma 4.13(a)). Under appropriate cicumstances also circuit-friendliness can be perpetuated from $G$ to $G'$ (Lemma 4.13(b)).

{\bf Lemma 4.13:} Let $(J',\Lambda')$ be the result of adding the path $[p_1,...,p_n]$ to $(J,\Lambda)$.	Furthermore, let  $G=(V,E)$ line-pres rank-model $(J,\Lambda)$ via  $J\ra E:r\mapsto\ol{r}$. Suppose that $p_1, p_n\in V$ are such that $\ol{p_1}$ and $\ol{p_n}$ are incident edges. 
\begin{enumerate}		
		\item[{(a)}]If $n\ge 2$ then there is a graph $G'$ line-pres rank-modeling $(J',\Lambda')$.
			\item[{(b)}] If $n\ge 3$ and $\ol{p_1},\ \ol{p_n}$ are such that $\Lambda$ has a line of type $\{p_1,p_n,q\}$ then, provided $G$ did so, $G'$ continues to rank-model circuit-friendly.
			\end{enumerate}

		Adopting graph theory parlance we say that a PLS has {\it small girth} if all its cycles consist of three or four lines. 

		{\bf Lemma 4.14:}	Let $(J,\Lambda)$ be a PLS.
		\begin{enumerate}
		\item[{(a)}] If $(J,\Lambda)$ is of small girth and has a line-pres rank-modeling graph, then $(J,\Lambda)$ is a BMPL.
		\item[{(b)}] If $(J,\Lambda)$ is an augmented UMP of small girth, then $(J,\Lambda)$ has a line-pres rank-modeling graph. 
		\item[{(c)}] If $(J,\Lambda)$ is an augmented UMP of type 1, then $(J,\Lambda)$ has a circuit-friendly rank-modeling graph. 
		\end{enumerate}

Parts (a),(b) can be rephrased as follows. When does a small girth PLS have a line-pres rank-modeling graph? It is necessary to be a BMPL. And sufficient to be an augmented UMP. Parts (b),(c) can be recast as follows. Do augmented UMPes have line-pres rank-modeling graphs?
Yes, if we add the small girth ingredient. Alternatively, we can add the type 1 ingredient and even get {\it circuit-friendly}  rank-modeling graphs. We stress that adding small girth or type 1 is merely sufficient for line-pres rank-modeling. If 'lucky'  one can do away with both of them and still get a line-pres rank-modeling graph. As a case in point, the augmented UMP $(J_2,\Lambda_2)$ is not of type 1 and has a large cycle (Figure 6b), yet {\it possesses} a line-pres rank-modeling graph $G_2$  (Fig. 6c). Unfortunately  $G_2$ is not circuit-friendly since the chordless circuit $\Gamma=\{9,11,12,14\}$  does not match the set of midpoints of a cycle in $(J_2,\Lambda_2)$. See also 8D.

					\section{ Matroids and modular lattices}

		 Both tight $k$-linear representations  of modular lattices, and tight embeddings of them into partition lattices (Sections 1, 2), fit the common hat of tight embeddings into lattices consisting of the closed subsets of matroids. Accordingly we investigate various kinds of matroids {\it modeling modular lattices}. This is akin to  Section 4 where we considered matroids {\it modeling PLSes}. But PLSes are absent in Section 5 and only return in Section 6.

	{\bf 5A Tight embedding of modular lattices into geometric lattices.} For any field $k$ let $LM(k^n)$ be the (modular) subspace lattice of $k^n$. A $k$-{\it linear representation} of $L$ is any homomorphism $L\rightarrow LM(k^n)$. Article [5] classifies up to isomorphism all $k$-linear representations of modular lattices which are acyclic in the sense of 6A. A bit more on that follows in Example 5.5 but by and large the {\it existence} of $k$-linear representations is  all we care about. Actually, generalizing $LM(k^n)$ we mainly look at $LM(K)$, which by definition is the lattice of closed sets of the matroid $M(K)$. Such lattices are also known as {\it geometric} lattices, see [10]. Whenever the matroid closure operator $\mathcal{P}(K)\ra\mathcal{P}(K):\ X\mapsto\widetilde{X}$ needs to be emphasized we write $M(K,\widetilde\ )$ rather than $M(K)$. In the special case where $K=k^n$ and $\widetilde{X}$ is the subspace generated by $X\subseteq k^n$, it can be notationally better to write $\langle X\rangle$ instead of $\widetilde{X}$. The 'Key-Lemma' below  is Lemma 5 in [15].

					{\bf  Lemma 5.1:} Let $L$ be a modular lattice with $J=J(L)$, and let $M(K,\widetilde{}\ )$ be a matroid. There is a tight embedding $\Phi:L\ra LM(K)$ (thus $d(LM(K)\ge d(L)$ allowed) iff the following holds. There is an injection $\varphi:J\ra K$ such that the induced submatroid $M(\varphi(J),-)$ is simple and such that (5.1) and (5.2) hold: 
					
					(5.1)\qquad $\ol{\varphi(J(a))}=\varphi(J(a))$ for all $a\in J$ {\it (dependency condition of the second kind)}.
					
					(5.2)\qquad $mrk(\varphi(J))=d(L)$ (\hbox{\it rank condition of the second kind})

This is good and well, but {\it how} does $\Phi$ arise from $\varphi$, and vice versa? Given $\varphi$ with (5.1) and (5.2), one can put $\Phi(a):=\widetilde{\varphi(J(a))}$. Conversely, given any tight homomorphism $\Phi$, for each $p\in J$ pick any $p'\in\Phi(p)\setminus\Phi(p_*)$ and define $\varphi:J\ra K$ by $\varphi(p):=p'$. Here $p_*\prec p$ is the unique lower cover of $p$ in $L$. It is crucial to distinguish the global closure $\widetilde{X}\ (X\subseteq K)$ from the induced closure $\ol{X}:=\widetilde{X}\cap\varphi(J)\ (X\subseteq \varphi(J))$. Thus $\ol{\varphi(J(a))}=\varphi(J(a))$ by (5.1) whereas generally $\widetilde{\varphi(J(a))}\neq\varphi(J(a))$!

{\bf Definition 5.2.} Let $M(E,-)$ be a simple  matroid and $L$ a modular lattice. We say that a bijection $\varphi: J(L)\ra E$ {\it lattice-models} $L$ if it satisfies (5.1) and (5.2). When we say that $M(E,-)$ {\it lattice-models} $L$, it is implied that such a $\varphi$ exists. $\Box$

As will be further investigated, condition (5.1) is the sibling of (4.1),
 and (5.2) the sibling of (4.2). Consider the three increasingly special cases where $M(E,-)$ incorporates  linear dependency (thus $E\subseteq k^m$ for some field $k$), or where particularly $k=GF(2)$, or where $M(E,-)$ is graphic. We then speak of $k$-{\it linear}, {\it binary} and {\it graphic} matroids modeling $L$. In the graphic case we usually speak (akin to Section 4) of  the {\it graph} $G=(V,E)$ {\it modeling} $L$.

{\bf 5B Trimming the Key-Lemma.} We first trim Lemma 5.1 to the $k$-linear case (Corollary 5.3, which is of lesser importance) and then, crucially, to graphs  in Theorem 5.4.

{\bf Corollary 5.3:}  Suppose the modular lattice $L$ has height $n=d(L)$ and is lattice-modeled by the $k$-linear matroid $M(E,-)$ via $\varphi:J(L)\ra E$. Then 
$\Phi(a):=\langle\,\varphi(J(a))\,\rangle\ (a\in L)$ provides a tight embedding $\Phi:L\ra LM(k^n)$. Conversely {\it each} tight embedding $\Phi':L\ra LM(k^m)$ comes from such a $\varphi$-induced tight embedding $\Phi:L\ra LM(k^n)$.

{\it Proof.} When $L$ is lattice-modeled by the $k$-linear matroid $M(E,-)$ one has $mrk(E)=d(L)=n$, and so $E$ can be viewed as a spanning submatroid of the matroid $K=k^n$. Lemma 5.1 hence yields a $k$-linear representation $\Phi:L\ra LM(k^n)$. Conversely, consider a map $\Phi'$ as above. Since  the subspace $\Phi'(L)$ of $k^m$ is isomorphic to $k^n$, we may view $\Phi'$ as a map  $L\ra LM(k^n)$. Applying Lemma 5.1 in the other direction shows that $\Phi'$ must be induced by a map $\varphi:J(L)\ra E$ satisfying (5.1) and (5.2). $\Box$

 {\bf 5B1 Partition lattices at last.} Recall from Section 1A that $Part(m)$ is  the lattice of all set partitions of $[8]$. The complete graph 
$CG(m)=([8],K)$ has vertex set $[8]$ and $|K|={m\choose 2}$ many edges. It is well known that the graphic  matroid $M(K,\widetilde{}\ )$ has rank $m-1$, and its height $m-1$ flat lattice 
$LM(K)$ is isomorphic to $Part(m)$ via $B\mapsto comp(B)$. Here for any edge set $B\subseteq K$ we define $comp(B)$ as the partition of $V$ whose parts are the vertex sets of the connected components of the subgraph $([8],B)$ of $CG(m)$. 

{\bf Theorem 5.4:} Let $L$ be a modular lattice of height $n=d(L)$. If the graph $G=([n+1],E)$ lattice-models $L$ via $\varphi: J(L)\ra E$ then $\Phi(a):=comp(\varphi(J(a))\ (a\in L)$ provides a tight embedding $\Phi:L\ra Part(n+1)$. Conversely, {\it each} tight embedding $\Phi':L\ra Part(m)$ has $m\ge n+1$ and 'comes from' such a $\varphi$-induced tight embedding $\Phi:L\ra Part(n+1)$.

The meaning of 'comes from' will be defined in a minute.

{\it Proof.}  That $\varphi$ as described induces a tight embedding  $\Phi:L\ra Part(n+1)$ follows from Lemma 5.1 and the remarks above. The converse direction is slightly more subtle than the converse direction  in Corollary 5.3. Thus let
$\Phi':L\ra Part(m)$ be  any tight embedding.
By Lemma 5.1 and the remarks above $\Phi'$ induces a graph $H'=([8],E)$ that lattice-models $L$. By gluing together vertices of potential disconnected components of $H'$ one gets a {\it connected} graph $H$. A moment's thought  shows that $H$ still lattice-models $L$. By connectedness $H$ now has $mrk(E)+1{(5.2)\atop =}d(L)+1=n+1$  vertices. Applying Lemma 5.1 in the other direction yields a tight embedding $\Phi:L\ra Part(n+1)$.  It is fair to say that $\Phi'$ {\it comes from} $\Phi$. $\Box$

{\bf Example 5.5.} Consider $L_2$ in Figure 21a. Recall its tight embedding $\Phi_2': L_2\ra Part(5)$  defined in Figure 2f of Section 2 where $5>d(L_2)+1$. As argued in the proof of Theorem 5.4 the tight embedding $\Phi_2'$ must be induced by  a lattice-modeling graph $H_2'$. One checks that $H_2'$ has two connected components, a triangle and a single edge (Figure 21b). The obtained connected modeling graph $H_2$ is shown in Figure 21c. Since it now has 4 vertices, it induces a tight embedding $\Phi_2:L_2\ra Part(4)$. 

\includegraphics[scale=0.45]{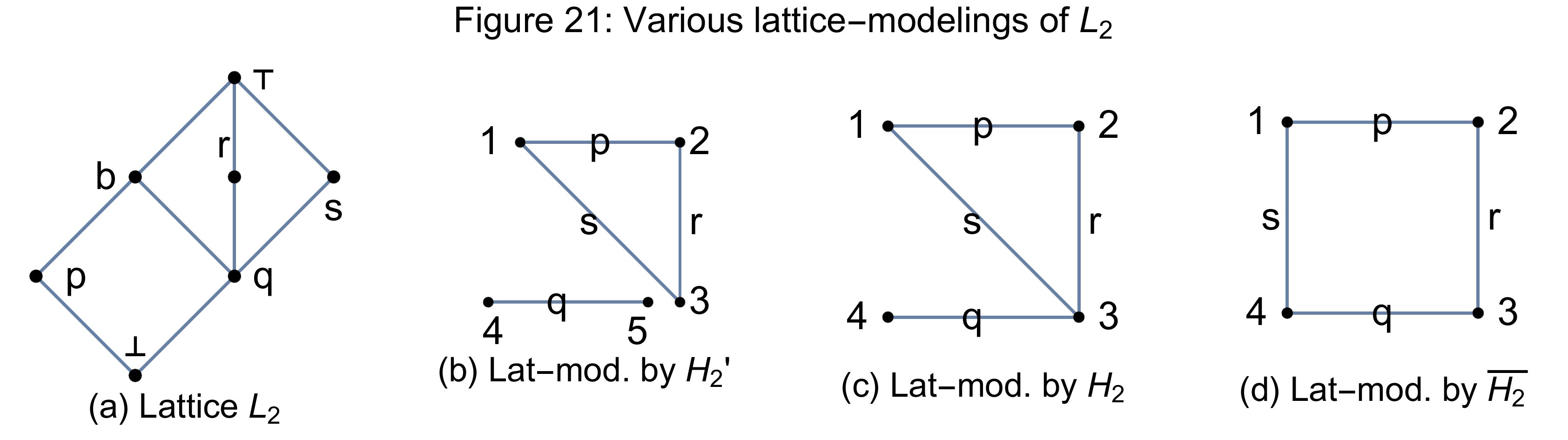}
				
We claim that $\overline{H_2}$ in Figure 21d is  another graph on four vertices that lattice-models $L_2$. Indeed, (5.2) is trivial and (5.1) will be verified as the variant $(5.1*)$ in a moment.
  Notice that $\overline{H_2}$ is 2-connected, whereas
 $H_2$ is not. Nevertheless,  both graphs yield the {\it same} embedding $\Phi_2$ of $L_2$. For instance in both cases $s\in L_2$ maps to the partition $\{\{2\},\{1,3,4\}\}$. This is  a surprise viewing that $s$ maps to different edges in 21c and 21d respectively. However, from the vantage point of $k$-linear (tight) representations $F_i: L_2\ra LM(k^3)$ the surprise is moderate. That is, by [5] any two $F_i,F_j$ must be {\it isomorphic} in the sense that for some vector space isomorphism $f:k^3\ra k^3$ it holds that $F_j(a)=f(F_i(a))$ for all $a\in L_2$. Whether for more complicated $L\neq L_2$ the isomorphism of all $k$-linear representations of $L$  still forces the isomorphism (in the natural sense) of all tight partition embeddings of $L$, remains an open question. When $L$ is not acyclic (6A) then its tight $k$-linear representations need no longer be isomorphic, and hence the classification of tight embeddings $L\ra Part(n+1)$ gets more complicated still. As previously stated, in the present article we are merely concerned with the existence of such embeddings (which is hard enough). $\Box$

{\bf 5B2. Reformulation of (5.1) in the graph case.} The criterion $(5.1*)$ below is a handy substitute for (5.1) in the graph case. Namely, let $G=(V,E)$ be a graph, $L$ a modular lattice, and $\varphi: J(L)\ra E$  a bijection. Then (5.1) is by [15, Lemma 10] equivalent to this condition:

$(5.1*)$\quad For each chordless circuit $\varphi(X)\subseteq E$ of $G$ it holds that $q\le \bigvee(X\setminus\{q\})$ for all $q\in X$.

For instance, the only chordless circuit $\{p,q,r,s\}$ of $\overline{H_2}$ in Figure 21d satisfies $(5.1*)$:

$p\le q\vee r\vee s,\ q\le p\vee r\vee s,\ r\le p\vee q\vee s,\ s\le p\vee q\vee r$

{\bf 5C Excursion to semimodular lattices.} A variation of Lemma 5.1 shows [15, Theorem 4] that even each {\it semi}modular lattice $L$ is tightly embeddable in a geometric lattice $LM(K)$, but possibly not in a 'nice' one like $LM(k^n)$ or $Part(n+1)$. That's because  the matroid  modeling $L$ is constructed as a submatroid of a certain matroid $(K,\widetilde\ )$ triggered by the
 {\it height function} $d:L\ra \N$. This makes $LM(K)$  quite unpredictable.
The proof of [15, Theorem 4] cuts short previous ones by Dilworth 1973 and Gr\"{a}tzer-Kiss 1986.

\section{Putting the pieces together and looking ahead}

As surveyed in 1B here we build our PLSes on the sets $J(L)$ of join-irreducibles of thin lattices $L$, and call them MoPLSes. Theorem 6.3 states that the existence of an augmented MoUMP of type 1 suffices for $L$ to be tightly partition embeddable.

{\bf 6A From PLSes to MoPLSes.} Let $L$ be a modular lattice. If $\ell\subseteq J(L)$  is maximal with respect to any distinct $p,q\in\ell$ yielding the {\it same} join $p\vee q$ (which we denote by $\overline{\ell}$), and if $|\ell|=3$ then $\ell$ is called a {\it line (of} $L$). (While $|\ell|>3$ occurs we are only interested, as for PLSes, in lattices $L$ all of whose lines have cardinality 3.)
Two lines $\ell$ and $\ell_0$ are {\it equivalent} if $\ol{\ell}=\ol{\ell_0}$. Any maximal family $\Lambda$ of mutually inequivalent lines yields a partial linear space $(J,\Lambda)$ since (3.1) is satisfied. 
We call it a {\it MoPLS}, where the Mo emphasizes that this PLS arises from a modular lattice. (In [5] the terminology 'base of lines of $L$' was used.)

{\bf Example 6.1.} Reconsider $L=L_1$ from Figure 1b, an upped version of which is Figure 22a. One has $J_1=J(L_1)=\{1',...,9'\}$. For instance $\ell=\{1',2',4'\}$ is a line with $\overline{\ell}=a_1$. The line $\{3',8',9'\}$, and six more, are equivalent to $\ell$. One verifies that the five lines in Figure 22b (=Figure 1c) yield a MoPLS $(J_1,\Lambda_1)$. $\Box$

\includegraphics[scale=0.52]{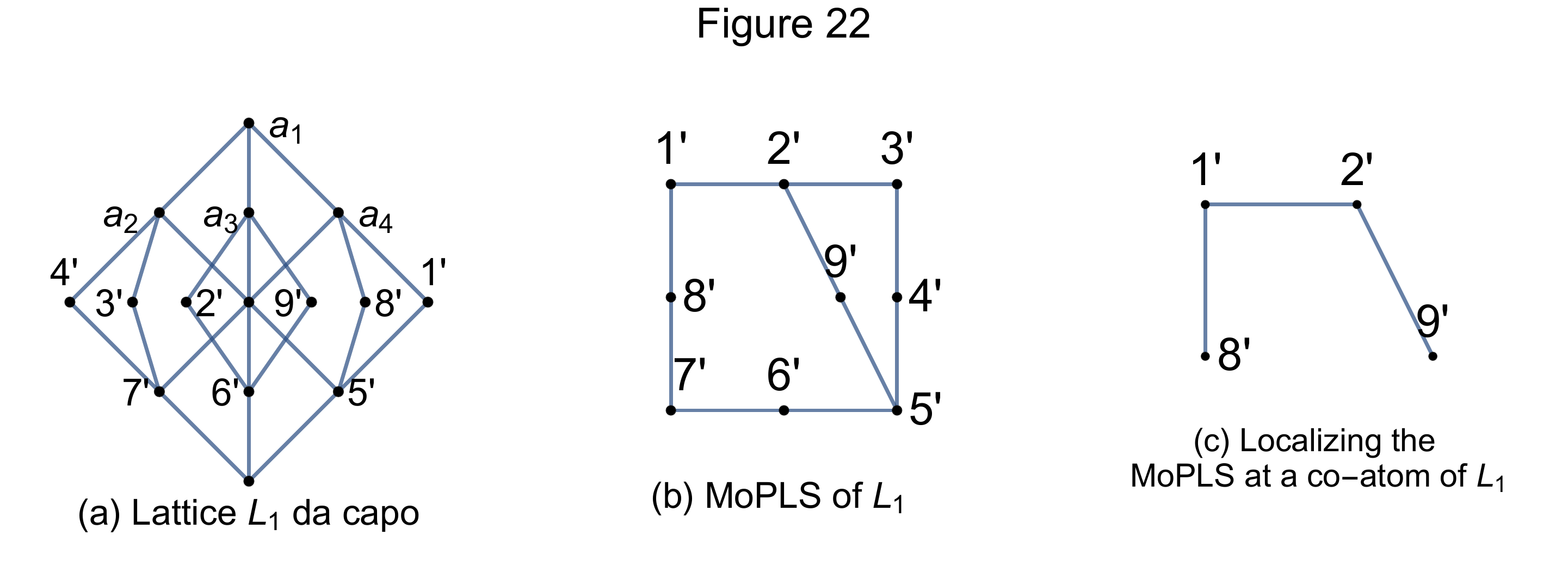}

  MoPLSes  extend Birkhoff's Theorem  in elegant ways to modular lattices $L$. Specifically, an order ideal $X$ of $(J,\le)$ is called $\Lambda$-{\it closed} if for all $\ell\in\Lambda$ it follows from $|\ell\cap X|\ge 2$ that $\ell\subseteq X$. If ${\CL}(J,\le,\Lambda)$ denotes the closure system of all $\Lambda$-closed order ideals then $a\mapsto J(a)$ turns out [5, Thm. 2.5] to be a lattice isomorphism from $L$ onto ${\CL}(J,\le,\Lambda)$. 
Each MoPLS of $L$ has $s(L)$ many connected components which are in bijection with the subdirectly irreducible factors of $L$. If $L$ is distributive then $s(L)=|J(L)|$, and so  ${\CL}(J,\le,\Lambda)$ boils down to $D(J,\le)$ in 2A.

A modular lattice is {\it acyclic} if all (equivalently: one) of its MoPLSes are acyclic. Inequality (2.1) is sharp exactly for acyclic modular lattices [5, 6.4]. See also Example 5.5. Akin to how 2-distributivity is a sweeping generalization of distributivity, local acyclicity (defined in 8G) is a sweeping generalization of acyclicity.

{\bf 6B Good news and bad news.}
In the sequel we focus on modular lattices which are  thin, i.e. (2D) they are  2-distributive and without covering sublattice $M_4$. According to [15, Lemma 19], and the conversion formula (8.1) in 8B, each MoPLS $(J,\Lambda)$ of a thin lattice $L$ satisfies

(6.1)\qquad $rk(J,\Lambda)=d(L)$.

For instance $L_1$ from Figure 22a is thin and satisfies $rk(J_1,\Lambda_1)=9-5=d(L_1)$. For modular lattices which are not 2-distributive equality (6.1) can fail, and probably always does. For instance $L=LM(GF(2)^3)$ has $(J_4,\Lambda_4)$ in Figure 11a as its unique MoPLS and (6.1) fails since $rk(J_4,\Lambda_4)=7-7<d(L)$.

  One may be led to say that Theorem 5.4  settles the tight embeddability of a modular height $n$  lattice $L$ into $Part(n+1)$: It works iff there is a connected graph $G=([n+1],E)$ and a bijection $\varphi: J(L)\ra E$ that obeys $(5.1*)$ and (5.2). But {\it how} can the existence of such a $\varphi$ be decided for concrete lattices? As argued in Section 2, for  $\varphi$ to exist it is necessary that $L$ be thin. As to sufficiency, that's why we turned to
	graphs and PLSes in Sections 3 and 4. Yet  graphs {\it line-pres modeling} PLSes merely constitute a crutch to the actually relevant, but more enigmatic graphs {\it lattice-modeling} thin lattices.
Let us look at this crutch more closely. 

The {\it good news} is that in view of link (6.1) the rank condition (4.2) for (Mo)PLSes is {\it equivalent} to the rank condition (5.2) for lattices. 

{\it Bad news} is that the dependency condition $(5.1*)$ for lattices is {\it not equivalent} to the dependency condition (4.1) for (Mo)PLSes. For instance $\overline{H_2}$ in Figure 21d satisfies $(5.1*)$ but we claim that (4.1) fails. Indeed, the only MoPLS of $L_2$ in Figure 21a consists of the line $\ell=\{p,r,s\}$ and the isolated point $q$. The line $\ell$ does not map onto a triangle of $\overline{H_2}$; in fact $\overline{H_2}$ features no triangles at all.
Hence there are bijections $J\ra E$ that lattice-model a thin lattice, but which do not line-pres (nor circuit-friendly) model any underlying MoPLS. The converse remains an open question:

{\bf Question 3}: Is there a 'freak' graph $G=(V,E)$ and a MoPLS $(J,\Lambda)$ of a thin lattice $L$ such that some bijection $\varphi: J\ra E$ line-pres rank-models $(J,\Lambda)$ yet does not lattice-model $L$?

The hoped for answer is 'no'. By e.g. investing the extra ingredient of circuit-friendliness we get what we want:

{\bf Theorem 6.2}: If the thin lattice $L$ has a MoPLS that is circuit-friendly rank-modeled by some graph, then $L$ is tightly partition embeddable.

{\it Proof.} Putting $J=J(L)$ and $G=(V,E)$ let $\varphi: J\ra E$ be a circuit-friendly rank-modeling. (We do not postulate (4.1), see  8D.) We need to verify $(5.1*)$ and (5.2). As to  (5.2), recall that this is automatic whenever $G$ satisfies (4.2). As to $(5.1*)$, let $\varphi(X)\subseteq E$ be a chordless circuit of $G$. By definition of circuit-friendliness in 4C6 the set $X$ is either a line $\{p,q,r\}$ or the set of $C$-midpoints of a cycle $C$ in $(J,\Lambda)$. For $X=\{p,q,r\}$ condition $(5.1*)$ holds in view of $p\vee q=p\vee r=q\vee r$. If $X$ is a set of $C$-midpoints then $(5.1*)$ holds by [15, Lemma 20]. $\Box$

 The crucial Lemma 20 in [15]  involves so called 'cycles of $M_3$-elements' (see 8E), thus a third type aside from cycles in PLSes and circuits in graphs. From Theorem 6.2 and Lemma 4.14(c) follows:

{\bf Theorem 6.3:} If the height $n$ thin lattice $L$ has an augmented  MoUMP of type 1 then there is a tight embedding $\Phi: L\ra Part(n+1)$.

{\bf 6C Graph-triggers.} Here we propose an alternative line of future research that omits messy cycles of $M_3$-elements, and is directed more towards matroid theory. 

{\bf Definition 6.4.} Let $M(E)$ be a matroid such that $E=E_1\cup E_2$ and $E=E_1\cap E_2=\{z\}$. Suppose each circuit $C$ of $M(E)$ is a circuit in either one of the submatroids $M(E_1)$ and $M(E_2)$, or it is of type

(6.2)\quad $C=(C_1\cup C_2)\setminus\{z\}$ where $C_i$ is a circuit of $M(E_i)$ and $C_1\cap C_2=\{z\}.$

In this situation $M(E)$ is isomorphic to the so-called {\it parallel connection} [10,p.240] of $M(E_1)$ and $M(E_2)$. $\Box$ 

One can show that when both matroids $M(E_i)$ are graphic then so is $M(E)$.

{\bf Definition 6.5.} Call the PLS $(J,\Lambda)$  a {\it graph-trigger} if each binary matroid $M(E)$ that line-pres rank-models $(J,\Lambda)$ must be graphic. $\Box$

The boring scenario that $(J,\Lambda)$ is a graph-trigger simply because  it cannot be line-pres rank-modeled by {\it any} binary matroid, will not disturb us henceforth.

{\bf Theorem 6.6:} Each acyclic PLS is a graph-trigger.

{\it Proof.} It suffices to prove the claim for connected PLSes, and this we do by induction on the number of lines. Thus let the connected $(J,\Lambda)$ be line-pres rank-modeled by the binary matroid $M(E)$ and put $n=rk(J,\Lambda)=mrk(E)$. We write $p'$ for the image of $p\in J$ in $E$. If $|\Lambda|=1$, then $|E|=|J|=3$ and $mrk(E)=rk(J,\Lambda)=2$. Hence $M(E)$ is isomorphic to the graphic matroid induced by a triangle. 

Suppose now that $|\Lambda|>1$. By acyclicity there is a $\ell\in\Lambda$ that intersects each other line in at most one point $z$ (and intersection $\{z\}$ {\it does} occur by connectedness). Putting $\ell=\{x,y,z\},\ J_1=J\setminus\{x,y\},\ \Lambda_1=\Lambda\setminus\{\ell\},\ J_2=\ell,\ \Lambda_2=\{\ell\}$, consider
 the PLSes $(J_1,\Lambda_1)$ and $(J_2,\Lambda_2)$, as well as their corresponding submatroids $M(E_1)$ and $M(E_2)$ of $M(E)$. As above it follows that $M(E_2)$ is induced by a triangle with edges $x',\ y',\ z'$. To unravel the structure of $M(E_1)$ first observe that $M(E_1)$ keeps on line-pres modeling $(J_1,\Lambda_1)$. Each acyclic PLS being sparse it follows from (4.2) that $mrk(E_1)\le rk(J_1,\Lambda_1)=n-1$. Yet the inequality cannot be strict because of 
$n=rk(J,\Lambda)=mrk(E)$ and $mrk(\{x',y',z'\})-mrk(\{z'\})=1$ and the submodularity of $mrk$. 
Hence $mrk(E_1)=rk(J_1,\Lambda_1)$, and so $M(E_1)$ rank-models $(J_1,\Lambda_1)$. By induction $M(E_1)$ must be graphic. By the remarks above it now suffices to show that each circuit $C$ of $M(E)$ which neither lies in $E_1$ or $E_2$ is of type (6.2). Evidently the case distinction below covers all cases.

{\it Case 1:} $x'\in C$ but $y'\not\in C$ (or dually with $x',\ y'$ switched). Then $C\setminus\{x'\}\subseteq E_1$, and so $x'$ is in the closure of $E_1$. This leads to the contradiction $mrk(E_1)=mrk(E_1\cup\{x'\})=n$.
Hence Case 1 is impossible.

{\it Case 2:} $x',\ y'\in C$. If $z'$ was in $C$, then $\{x',y',z'\}$ was a dependent subset of $C$, which by definition of a circuit implies $C=\{x',y',z'\}$. This is impossible because $C\not\subseteq E_2$. To fix ideas, say $C=\{x',y',a',b',c'\}$ (and $z'\not\in C$). Since $M(E)$ is binary we may think of these elements as lying in some $GF(2)$-vector space. It holds that $x'+\cdots+c'=0$ but no proper subset of $C$ sums to $0$. Hence $z'+a'+b'+c'=(x'+y')+a'+b'+c'=0$ but no proper subset of $\{z',a',b',c'\}$ sums to $0$. It follows that $C_1=\{z',a',b',c'\}$ and $C_2=\{x',y',z'\}$ are circuits as required in (6.2). $\Box$

The relevance of graph-triggers derives from Theorem 6.7 below. Recall from Question 3 that we do not know whether each graphic matroid  that line-pres rank-models a MoPLS $(J(L),\Lambda)$ via $\varphi:J(L)\ra E$ automatically lattice-models $L$. If $\varphi$ fails, is there another $\varphi'$ that models both? Nobody knows. Interestingly, when we switch from graphic to binary matroids, the answer is yes. Specifically, let $M(E)$ be a binary matroid. Call a bijection $\varphi: J(L)\ra E$ {\it bi-modeling} if it simultaneously line-pres rank-models $(J(L),\Lambda)$ and lattice-models $L$.

{\bf Theorem 6.7:} Let $L$ be a thin lattice. If $L$ has a MoPLS which is a graph-trigger, then $L$ is tightly embeddable into a partition lattice.

{\it Proof.} Let $(J,\Lambda)$ be a MoPLS of $L$ which is a graph-trigger. As for any MoPLS of $L$, according to [5, Thm. 5.1] (see also 8G) there is a  binary matroid $M(E)$ and a bi-modeling bijection $\varphi: J\ra E$. Because $\varphi$ line-pres rank-models $(J,\Lambda)$, and $(J,\Lambda)$ is a graph-trigger, $M(E)$ is in fact graphic. Therefore, because $\varphi$  lattice-models $L$, Theorem 5.4 yields a tight embedding into $Part(n+1)$ for $n=d(L)$. $\Box$

Perhaps the proof of Theorem 6.6 can be adapted to affirmatively answer 
 
{\bf Question 4:} Is each QIMP a graph-trigger?

If yes, then in view of Lemma 3.2 chances are that UMPes are graph-triggers as well.
Emboldened by the partial success in Lemma 4.14 we therefore ask:

{\bf Question 5:} Is each (augmented) UMP, or even each BMPL a graph-trigger?

If yes, we would have an elegant proof of Theorem 6.3 which currently relies on the ugly $M_3$-cycles in Lemma 20 of [15]. When tackling Questions 4 and 5 the recent criterion [M,Thm.1] for establishing the graphicness of certain binary matroids $M(E)$ looks promising. What makes the Mighton-test appealing is that it suffices to target any {\it one} suitable base of $M(E)$. It is also possible that  instead of piling up benign midpoint-links, sparsity (3E)
provides a more fruitful inductive scheme to settle parts of Question 5 affirmatively.
 If even BMPLes are graph-triggers, then combinatorists can hand the flag to modular latticians. They must adapt the non-embeddability argument for $L_1$ (see 8E) to arbitrary modular lattices without MoBMPLes. If successful, the tight partition embeddability problem would finally be solved.

\section{The longer proofs}

We focus on the bare proofs of Theorem 2.1 and Lemmata 3.8, 4.7, 4.13, 4.14. But occasionally additional background information contained in Section 8 will be pointed out.

{\bf  Proof of Theorem 2.1.} Recall that the $\theta_i$'s bijectively match the perspectivity classes of prime quotients. Thus for each covering pair $a\prec b$ in $N$ there is {\it exactly one} $\theta_i$ with $(a,b)\not\in \theta_i$. 
 Let now $f_i:N_i\ra Part(n_i)$ be tight embeddings $(1\le i\le s)$. In order to get a tight embedding $f:N\ra Part(n)$ let 

$$T:N\ra N_1\times\cdots\times N_s,\quad a\mapsto (T_1(a),\ldots,T_s(a))$$

be a subdirect embedding of $N$. Thus $T$ is injective and all component maps $T_i:N\ra N_i$ are surjective. In order to see that the homomorphism 
$f(a):=(f_1(T_1(a)),\ldots,f_s(T_s(a))$ is a tight embedding of $N$ into $Part(n_1)\times\ldots\times Part(n_s)$, take any covering pair $a\prec b$ in $N$. By the injectivity of $T$ there is at least one index, say $i=1$, such that $T_1(a)<T_1(b)$. From $a\prec b$ and the surjectivity of $T_1$ follows that in fact $T_1(a)\prec T_1(b)$. As observed at the beginning of the proof, there is {\it no other} $j$ with $T_j(a)\prec T_j(b)$, thus $T_j(a)=T_j(b)$ for all $j>1$. Since $f_1$ is cover-preserving, one has $f_1(T_1(a))\prec f_1(T_1(b))$, and so
$f(b)=(f_1(T_1(b)),f_2(T_2(a)),\ldots,f_s(T_s(a))$ is indeed an upper cover of $f(a)$. Finally, putting $n:=n_1+\cdots +n_s$  observe that $Part(n_1)\times\ldots\times Part(n_s)$ is isomorphic to 
the sublattice of  $Part(n)$ that consists of all partitions refining $\{\{1,\ldots,n_1\},\{n_1+1,\ldots,n_1+n_2\},\ldots,\{..,n-1,n\}\}$. $\Box$

\vskip 1.6cm

	{\bf  Proof of Lemma 3.8.}  By induction it suffices to show the following. Let $(J',\Lambda')$ be the result of adding a benign type 1 midpoint-link $P'$ to some cycle $C_0$ of a  BMPL $(J,\Lambda)$. By way of contradiction assume that $(J',\Lambda')$ possesses a cycle with a {\it bad}, i.e. non-benign midpoint-link. We will show that this forces the existence of a bad midpoint-link already in $(J,\Lambda)$. In fact the bad midpoint-links in $(J,\Lambda)$ will be seen to have the same type (1 or 2) as the triggered one in $(J',\Lambda')$. In particular this will prove part (b) along with part (a). We now embark on part (a).
	
	Specifically 
		let	 $P'$ above be $P'=[x,z_1,...,z_s,y]$. The points in $P'^*\setminus\{x,y\}$ are the {\it inner} points of $P'$. Note that the possibility $s=0$ in 3F (i.e. $P'=[x,y]$)  now becomes $s\ge 1$ since a type 1 midpoint-link comprises at least two lines. Furthermore type 1 implies that $[x,y]=\{x,y,u\}$ is a line of $C_0$. We call it the line {\it supporting} $P'$. Further let $C'$ be a cycle of $(J',\Lambda')$ with a {\it bad} midpoint-link $Q'=[r_1,...,r_n]$. (In Figures 8a, 8b, 8c the $C'$ is the 'hexagon' and one has $n=3,5,4$ respectively.) Since by assumption $(J,\Lambda)$ is a BMPL, some lines of $C'$ or/and $Q'$ belong to $P'$. A line of $C'$ or $Q'$ that belongs to $P'$ is called {\it new}. Because $r_n$ in Figure 8a is not quasi-isolated, yet all inner points of $P'$ are quasi-isolated, it cannot be that all dashed lines in Figure 8a are new.  It follows that all new lines  belong to either $Q'$ (Case 1), or they all belong to $C'$ (Case 2).

			\includegraphics[scale=0.55]{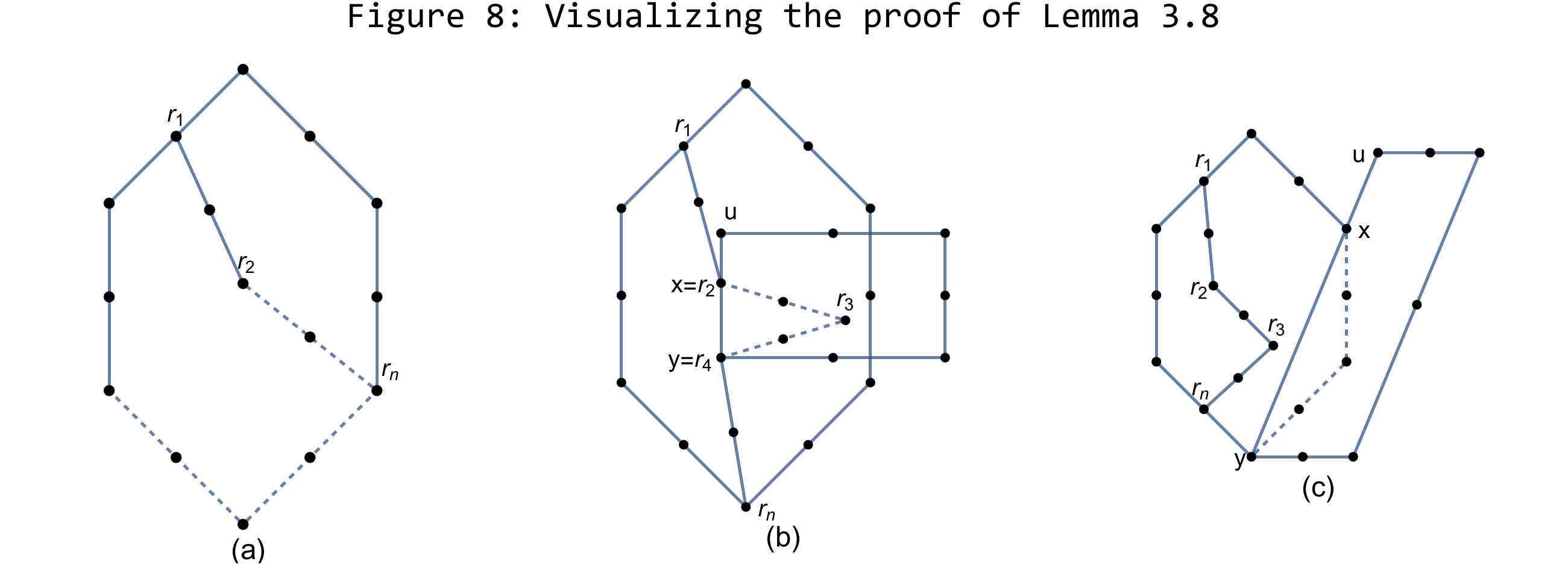}

	{\it Case 1:} All new lines  belong to $Q'$. Then again by the quasi-isolation of  inner points of $P'$ the {\it whole} of $P'$ must be part of $Q'$. See Figure 8b where $P'$ is shown dashed. Recall that $P'=[x,...,y]$ is a type 1 benign midpoint-link of some cycle $C_0$ of $(J,\Lambda)$ (which is the 'square' in Figure 8b and the 'parallelogram' in 8c).  Let us shrink the path $Q'=[r_1,r_2,r_3,r_4,r_5]$ to $Q=[r_1,r_2,r_4,r_5]$, where $[r_2,r_4]=[x,y]=\{x,y,u\}$ is the line of $C_0$ that supports $P'$. If we can argue that $[r_2,r_4]$ (or the likes) is not a $C'$-line then we get the desired contradiction that in $(J,\Lambda)$ the cycle $C:=C'$ has the bad midpoint-link $Q$. To show this, note that while $u\in C'^{*}$ is possible, we have $r_2,\ r_4\not\in C'^{*}$ (by definition of a midpoint-link $[r_1,...,r_n]$ of $C'$), and so $[r_2,r_4]$ is not a $C'$-line. If say $P'$ stretched from $x=r_2$ to $y=r_n$ then $[r_2,r_n]$ still is no $C'$-line since $r_2\not\in C'^{*}$. If  $P'$ stretched from $x=r_1$ to $y=r_n$ then $r_1,\ r_n\in C'^{*}$ but $u\not\in C'^{*}$. Thus again $[r_1,r_n]$ is no $C'$-line.

	{\it Case 2:} All new lines belong to $C'$. Again by the quasi-isolation of  inner points of $P'$ the {\it whole} of $P'$
	must be part of $C'$. We further branch according to the behaviour of the supporting line $[x,y]$ of $P'$.
	
	{\it Subcase 2.1:} The line $[x,y]$ is not a line of $C'$ or $Q'$, see Figure 8c where $P'$ is again rendered dashed. Consider the cycle $C$ obtained from $C'$ by shrinking $P'=[x,...,y]$ to the line $[x,y]$. Then {\it in} $(J,\Lambda)$ the cycle $C$ has the bad midpoint-link $Q:=Q'$. This still holds when (other than in Figure 8c) one has $\{x,y\}\cap\{r_1,r_n\}\not=\emptyset$. But$\{x,y\}=\{r_1,r_n\}$ is impossible since $x,\ y$ are both $C'$-junctions, whereas at least one of $r_1,\ r_n$ is a $C'$-midpoint.
					
	{\it Subcase 2.2:} The line $[x,y]$ is simultaneously a line of $Q'$. Since in $Q'=[r_1,...,r_n]$ only $r_1,\ r_n$ sit on $C'$, it follows that $\{x,y\}=\{r_1,r_n\}$, which yields the same contradiction as in Subcase 2.1.
	 Hence Subcase 2.2 is impossible.

	{\it Subcase 2.3:} The line $[x,y]$ is  a line of $C'$. Subcase 2.3.1: $[x,y]$ is a line of $P'=[x,z_1,...,z_s,y]$. This is impossible since $s\not =0$ as previously noticed. Subcase 2.3.2: $[x,y]$ is not a line of $P'$. Then the $C'$-lines are exactly $[x,y]$ and the lines of $P'$. Hence at least one of $r_1,\ r_n$ is an inner point of $P'$, which contradicts these points being quasi-isolated.	
	 $\Box$
					
		\vskip 1.6cm

{\bf  Proof of Lemma 4.7.} Suppose the cycle  $C$ in Lemma 4.7 is $C=(p_1,\ldots,p_n)$. For each $r\in C^*$ let $r':=\psi(r)$ be the associated edge.
	
	{\bf Claim A.}\qquad $\{p_1',\ldots,p_n'\}=star(v)$ for some vertex $v$ of $G$.
	
(Recall that $star(v)$ consists of all edges incident with $v$.) Since each $\{p_i,q_i,p_{i+1}\}\in\Lambda^*$  yields a triangle $\{p_i',q_i',p_{i+1}'\}$ of $G$, it will follow from Claim A that $\psi(C^*)$ is a wheel, thus proving Lemma 4.7. The proof of Claim A will be based upon:
	
	{\bf Claim B.}\qquad If Claim A fails then the edge set $\{p_1',\ldots,p_n'\}$  contains the edge set of a circuit of $G$.
	
	{\it Proof of Claim B.} 	Because for each line $\{p_i,q_i,p_{i+1}\}$ the edge set $\{p_i',q_i',p_{i+1}'\}$is a triangle of $G$, each edge $p_{i+1}'$ is incident with edge $p_i'$ (modulo $n$). Say $p_2'$ is incident with $p_1'=\{v_1,v_2\}$ in $v_2$. Since Claim A fails there is $i\ge 2$ such that $p_i'$ is incident with $v_2$ but $p_{i+1}'$ is not, see Figure 17. If $p_{i+1}'$ is incident with $v_1$ then $\{p_1',p_i',p_{i+1}'\}$ is a triangle (whence circuit) of $G$. Otherwise consider $p_{i+2}'$. As shown in Figure 17 there are two options for $p_{i+2}'$. If $p_{i+2}'$ is incident with $p_1'$ then we get again a circuit of $G$; if not continue with $p_{i+3}'$, and so on. Because at the latest $p_n'$ is incident with $p_1'$, there must be a cycle in $G$. This proves Claim B.

	\includegraphics[scale=0.57]{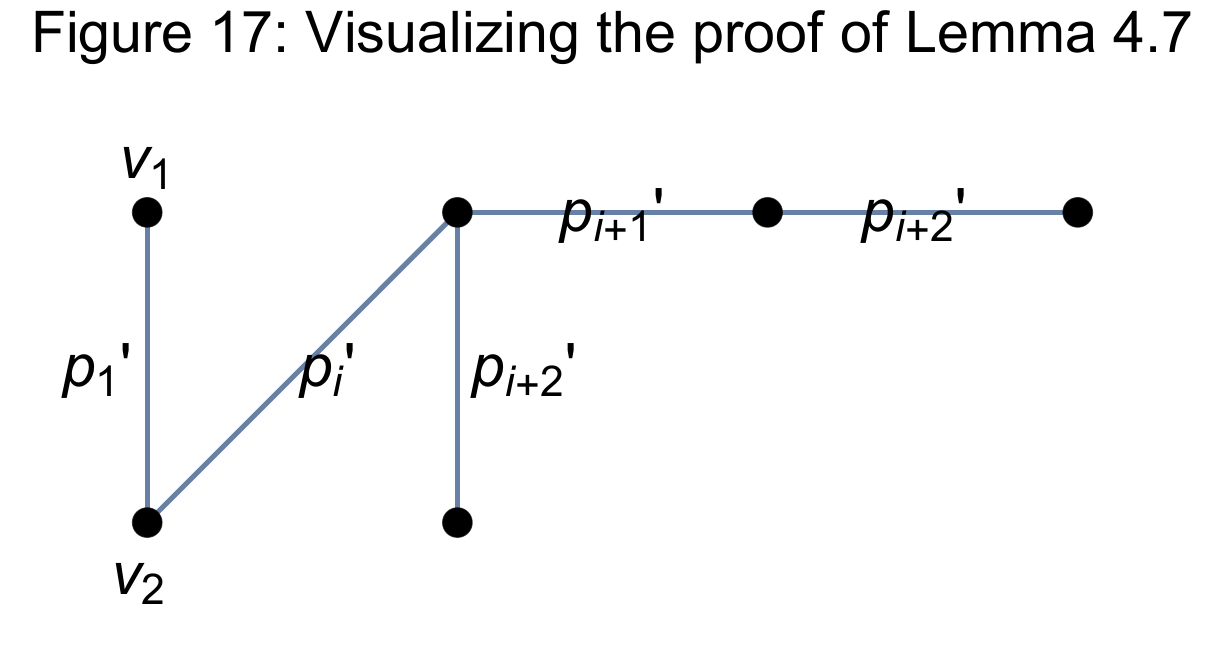}

{\it Proof of Claim A.} By way of contradiction assume that Claim A fails. On the one hand $\{p_1',\ldots,p_n'\}$ spans the submatroid of $M(E)$ induced by $\psi(C^*)$ since each $q_i'$ is in a circuit $\{p_i',q_i',p_{i+1}'\}$. On the other hand $\{p_1',\ldots,p_n'\}$ is dependent since it contains a circuit by Claim B. Hence $mrk(\psi(C^*))=mrk(\{p_1',\ldots,p_n'\})\le n-1$.
 Yet $rk(C^*,\Lambda^*)=|\Lambda^*|=n$, and so $rk(C^*,\Lambda^*)\neq mrk(\psi(C^*))$. This contradicts the assumption in Lemma 4.7, and thus proves Claim A. $\Box$

	\vskip 1.6cm
	
	{\bf  Proof of Lemma 4.13.} {\bf (a)} The case $n=2$ being dealt with in Subsection 4C7 we focus on $n\ge 3$. Let $0\in V$ be the vertex in which $\ol{p_1}$ and $\ol{p_n}$ intersect. Add {\it new} spokes $\ol{p_2},...,\ol{p_{n-1}}$ centered in $0$ (recall $n-1\ge 2$) and add rims to get a  'dented wheel'. (This means the following. Unless there happens to be an edge $\ol{q}\in E$ such that $\{\ol{p_1},\ol{p_n},\ol{q}\}$ is a triangle in $G$ (as e.g. in (b)), our wheel is lacking exactly the one rim between the spokes $\ol{p_1}$ and $\ol{p_n}$.) The new edges yield a blown up graph $G'=(V',E')$. For instance in Figure 19a we have $V'=V\cup\{\beta_2,...,\beta_{n-1}\}$ and $E'\setminus E$ consists of the boldface (non-dashed) edges. Generally $G'$ line-pres models $(J',\Lambda')$ since each $\{p_i,q_i,p_{i+1}\}\in\Lambda'\setminus\Lambda$ maps to a triangle in $G'$. It remains to check that the PLS-rank and matroid-rank increase the same amount, i.e. that
			
			$\rho:=rk(J',\Lambda')-rk(J,\Lambda)=mrk(E')-mrk(E)=:\mu.$
			
			Indeed, $\mu=|V'|-|V|=n-2$ matches $\rho=|J'|-|\Lambda'|-(|J|-|\Lambda|)=(|J'|-|J|)-(|\Lambda'|-|\Lambda|)=(2n-3)-(n-1)=n-2$ (recall $p_1,\ p_n\in J$).
		
		\includegraphics[scale=0.73]{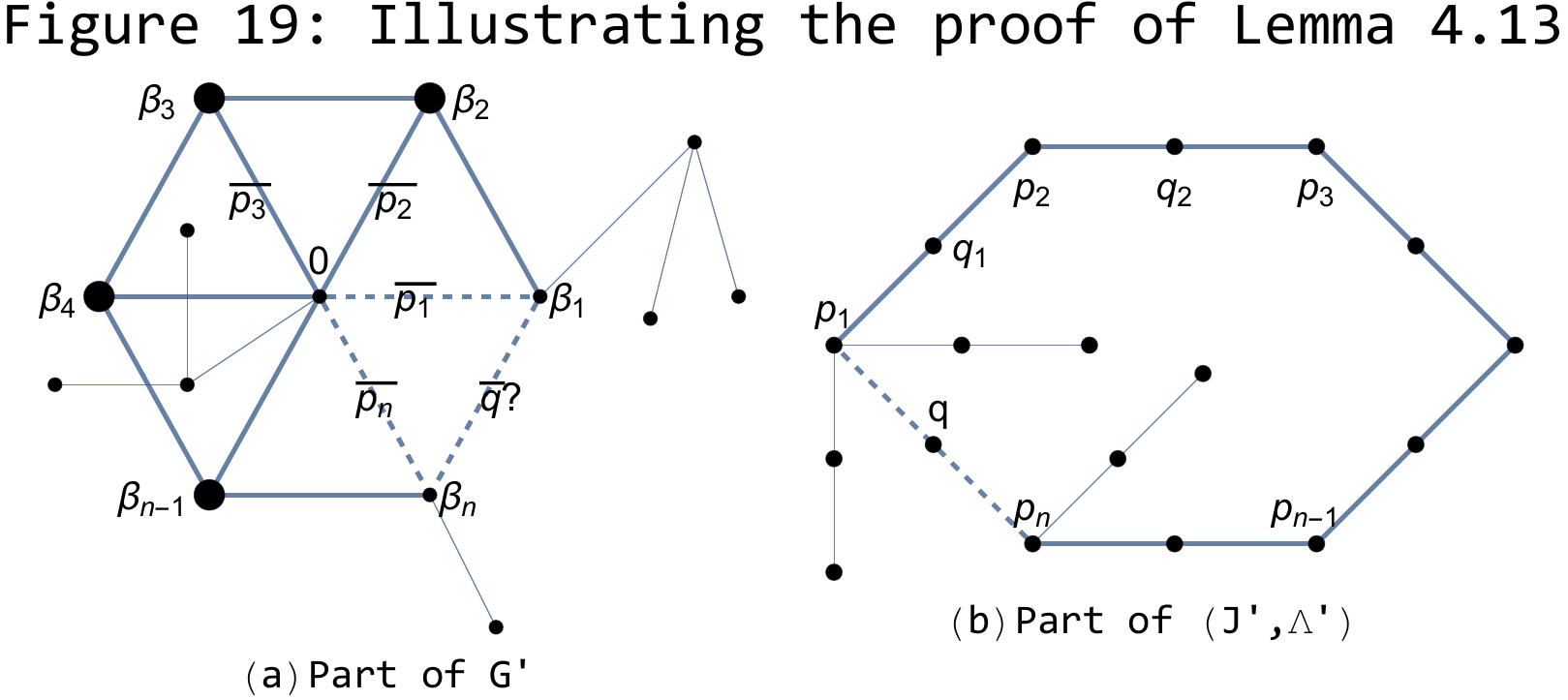}	
			
	{\bf (b)} By assumption $\{p_1,p_n,q\}\in\Lambda$ (indicated dashed in Figure 19b). This line maps to the dashed triangle in Figure 19a. Let now $\Gamma$ be any chordless circuit in $G'$. We need to show that $\Gamma$ either matches a line or the midpoints of a cycle in $(J',\Lambda')$. Since $G$ models $(J,\Lambda)$, we can assume that $\Gamma\not\subseteq E$.
			
			{\it Case 1:} There is a new spoke $\ol{p_j}$ in $\Gamma$, i.e. $1<j<n$. Hence $\ol{p_j}=\{0,\beta\}$ for some $\beta\in V'\setminus V$; such as $\beta=\beta_3$ in Fig 19a. There are only two possibilities for the  neighbour $\not =0$  of $\beta_3$ within $\Gamma$, namely $\beta_2$ or $\beta_4$. Thus $\Gamma$ is either of type $\Gamma=\{\{0,\beta_3\},\{\beta_3,\beta_2\},\ldots\}$ or $\Gamma=\{\{0,\beta_3\},\{\beta_3,\beta_4\},\ldots\}$.
			Being chordless forces $\Gamma$ in both cases to be a triangle (thus $\{\{0,\beta_3\},\{\beta_3,\beta_4\},\{\beta_4,0\}\}$ in the second case), and in both cases this triangle matches a line in $(J',\Lambda')$.
			
			{\it Case 2:} All edges in $\Gamma\setminus E$ (and there {\it are} such) are rims of the constructed wheel. To fix ideas, say $\{\beta_3,\beta_4\}\in\Gamma$ in Figure 19a. 
			Since $\Gamma$ contains no spokes, it follows that $\{\beta_i,\ \beta_{i+1}\}\in\Gamma$ for all $1\le i\le n-1$. {\it Subcase 1:} $\ol{q}\in\Gamma$. Then $\Gamma=\{\{\beta_1,\beta_2\},\ldots,\{\beta_n,\beta_1\}\}$, and so $\Gamma$ matches the set of midpoints $\{q_1,\ldots,q_n,q\}$ of a cycle in $(J',\Lambda')$, see Fig. 19b. {\it Subcase 2:} $\ol{q}\not\in\Gamma$. Then $\Gamma$ would have $\ol{q}$ as a chord, which is impossible. $\Box$

\vskip 1.6cm

{\bf  Proof of Lemma 4.14.} {\bf (a)} By Lemma 4.11 it suffices to show this:

{\bf Claim C.} If $W=(V,E)$ line-pres models a cycle-PLS $(C^*,\Lambda^*)$ with $|\Lambda^*|\le 4$ then $W$ must be a wheel.

	{\it Proof of Claim C.} Suppose first (case A) that strict inequality $<$ takes place in (4.2). Then $|V|-1=mrk(E)<rk(C^*,\Lambda^*)=|\Lambda^*|$. If $|\Lambda^*|=3$ then $|V|\le 3$, whence $|E|\le {3\choose 2}<6=|C^*|$. If $|\Lambda^*|=4$ then $|V|\le 4$, whence $|E|\le {4\choose 2}<8=|C^*|$. In both subcases this contradicts 
	the bijectivity of $\psi: C^*\ra E$. (This cardinality argument breaks down when  $|\Lambda^*|=5$ since then $|E|\le {5\choose 2}=10=|C^*|$. And indeed things can go wrong, as witnessed by Figure 16.) Now suppose (case B) that equality takes place in (4.2). Then the claim follows from Lemma 4.7. 
	
{\bf (b)} The claimed modeling graph  is not uniquely determined (not even up to isomorphism). It depends on the particular way we piled up, starting with a UMP, our benign midpoint-links. In any case, by induction it suffices to show the following. Let $(J,\Lambda)$ be a PLS of small girth that is line-pres rank-modeled by $G=(V,E)$ via $\psi: J\ra E$. If $(J',\Lambda')$ is obtained by adding a benign midpoint-link $P$ to a cycle $C$ of $(J,\Lambda)$ then there is some graph $G'$ modeling $(J',\Lambda')$.

{\it Case 1:} $P=[p,...,q]$ is of type 1. Thus $p$ and $q$ are on a common line $\ell$ of $C$, being the $C$-junction and $C$-midpoint respectively. Since $\psi(\ell)$ is a triangle in $G$, the edges $\ol{p}=\psi(p)$ and $\ol{q}=\psi(q)$ are incident. It follows from Lemma 4.13(a) that some graph $G'$ models $(J',\Lambda')$.

{\it Case 2:} $P=[q_1,...,q_2]$ is of type 2. Then there are lines $\ell_1$ and $\ell_2$ of $C$ such that $q_1$ and $q_2$ are the $C$-midpoints of
 $\ell_1$ and $\ell_2$ respectively, see Figure 20a.

	\includegraphics[scale=0.73]{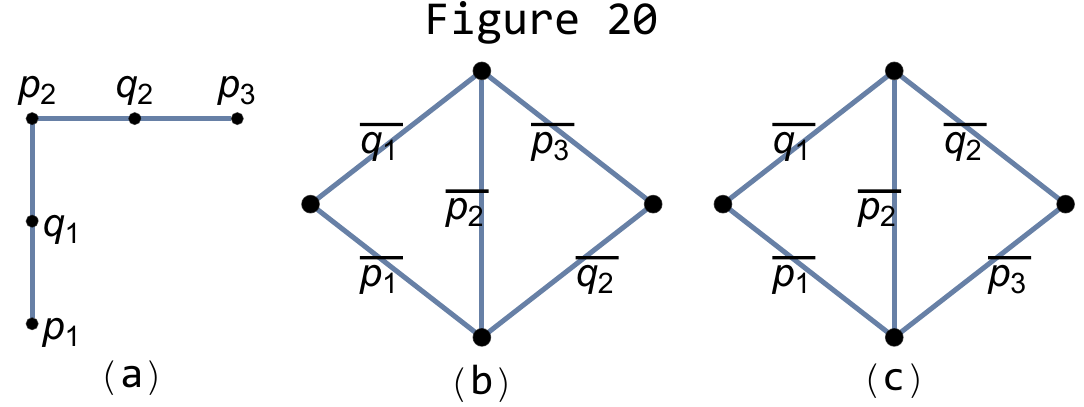}
	
	The associated edges in $G$ either behave as in Figure 20b or 20c. In order to again invoke Lemma 4.13(a) we must make sure that 20c takes place, because only then $\ol{q_1}$ and $\ol{q_2}$ are incident. (This is akin to flipping triangles in 4C1.) That's where we need that the induced cycle-PLS $(C^*,\Lambda^*)$ has $|\Lambda^*|\le 4$. It then follows as in the proof of (a) that $\psi(C^{*})$ is a wheel, and so $\ol{q_1}$ and $\ol{q_2}$ are incident rims.
	
	{\bf (c)} The proof is as in Case 1 of (b), with the addition that by induction (rooted in the starter UMP, see Lemma 4.12) we can assume that $G$ models $(J,\Lambda)$ circuit-friendly. By Lemma 4.13(b) this carries over to $G'$ modeling $(J',\Lambda')$.  $\Box$

\section{Further remarks}

{\bf 8A A garden of trees.} In 4C2 it was claimed that 'we may assume...'. The justification is as follows. If a QIMP is not 2-connected (as Figure 3c) then it is a tree  of its 2-connected components in the sense of 3C1. Thus a tree of QIMPes (2-connected or not) need not be a proper UMP, it can remain a QIMP. This e.g. happens if always {\it junctions} of component QIMPes are glued with {\it junctions} of component QIMPes (as opposed to midpoints or no name elements). While each tree of UMPes evidently stays a UMP, it is less clear whether a tree of augmented UMPes is again an an augmented UMP. In any case, it is a BMPL because generally trees of BMPLes are again BMPLes.

{\bf 8B About the PLS-rank.} In $(J_1,\Lambda_1)$ of Figure 1c imagine splitting the point $1'$ in half; one half maintains the line $\{1',2',3'\}$, the other half maintains the line $\{1',7',8'\}$. In the resulting ten points PLS split $2'$ in similar fashion. This yields an eleven points PLS which is acyclic. Alternatively, start again with $(J_1,\Lambda_1)$ and split $5'$ in order to 'detach' the line $\{5',6',7'\}$. In the resulting PLS split $5'$ once more in order to detach the line $\{2',5',9'\}$. This yields another acyclic PLS on eleven points. The precise definition of  {\it point splitting} is given in [15,p.216]. There it is further shown that for each PLS $(J,\Lambda)$, if not the point splittings themselves, at least the {\it number} of point splittings required to turn $(J,\Lambda)$ into an acyclic PLS $(J',\Lambda')$ with $c(J',\Lambda')=c(J,\Lambda)$ is an invariant $r^*(J,\Lambda)$ of $(J,\Lambda)$. 
In [15] the rank of a PLS was defined as $r(J,\Lambda):=|\Lambda|-r^*(J,\Lambda)$. Do not confuse this with $rk(J,\Lambda)$ as defined in (3.3). Pleasantly, as
noticed by Jim Geelen, for {\it our} PLSes with exclusively 3-element lines one can get rid of $r(J,\Lambda)$. Specifically we claim that

 (8.1) \qquad $rk(J,\Lambda)=|\Lambda|-r^*(J,\Lambda)+c(J,\Lambda)\ \ (=r(J,\Lambda)+c(J,\Lambda)).$

In order to prove (8.1), note that both sides in (8.1) agree for one-point PLSes in view of $1-0=0-0+1$. Since both sides of (8.1) are clearly  additive for disconnected  PLSes, it suffices to show that besides rk also the right hand side of (8.1) satisfies $(3.4),(3.5),(3.6)$. Thus let $x$ be the value of the right hand side for a connected PLS having 
point set $T_i:=\ell_1\cup\ldots\cup\ell_i$. Let us check that upon adding $\ell_{i+1}$ the right hand side changes as it must. {\it Case 1:} $|T_i\cap \ell_{i+1}|=1$. Then $x$ becomes $x+(1-0+0)=x+1$. {\it Case 2:} $|T_i\cap \ell_{i+1}|=2$. Then $x$ becomes $x+(1-1+0)=x$. {\it Case 3:} $|T_i\cap \ell_{i+1}|=3$. Then $x$ becomes $x+(1-2+0)=x-1$.

Note that the 'clumsy' rank $r(J,\Lambda)$ cannot be avoided for PLSes with larger lines. For instance, if $(J,\Lambda)=(J,\{\ell\})$ with $|\ell|=7$ then  $r(J,\Lambda)+c(J,\Lambda)$ correctly yields $(1-0)+1=2$ whereas  $rk(J,\Lambda)=7-1=6$ is way off.

{\bf 8C The 1992 ancestors of QIMPes, UMPes and BMPLes.} In [15] QIMPes are called 'mpi' (=middle point isolated), and UMPes are called 'regular'. In [15, Thm.23] regular PLSes were generalized to 'quasiregular' PLSes in an awkward manner. In contrast Theorem 6.3 of the present article generalizes UMPes more neatly  to augmented (Mo)UMPes of type 1.
While haphazard midpoint-links appear already in [15], credit goes to Manoel Lemos for shining more light on that matter in Lemma 4.11 which inspired the author: Upon reckognizing that cycle-preservation cannot be taken for granted in Lemos' original argument, one thing led to another.  For picture aficionados, page 226 of [15] showcases a 71-element thin lattice $L$ of height 11 and its
 quasiregular MoPLS, and  [15,p.241] explicitely lists  how the elements of $L$ embed into $Part(12)$. Up to an uninteresting line this MoPLS coincides with the augmented UMP of type 1 in Figure 4c ({\it with} dashed line).

	{\bf 8D Further deliberations on circuit-friendliness.} Recall that line-preservation (4.1) is not required in the statement of Theorem 6.2. Let us discuss the pros and cons of adding (4.1). If (4.1) is absent then we have more freedom to design $G$ and a circuit-friendly model $\psi:J\ra E$. On the other hand, freedom is a curse if one doesn't know how to exploit it. Put another way, the presence of (4.1) may be beneficial as a guidline to construct a circuit-friendly $\psi$, if only by means of the naive algorithm in 4C1. We can be a bit smarter though. To begin with, recall that by definition in a chordal graph all chordless circuits are triangles. Generalizing 'chordal', call a graph {\it wheely} if all chordless circuits $C$ derive from wheels, i.e $C$ is either a triangle or the rim-set of a nondegenerate wheel.

	{\bf Question 6:} What structure theorem can be proven about wheely graphs?

	This  is relevant because when $G=(V,E)$ circuit-friendly {\it and} line-preserving models $(J,\Lambda)$ via $\psi:J\ra E$ then $G$ must be wheely. In order to prove it, let $\Gamma\subseteq E$ be a chordless circuit. If $\psi^{-1}(\Gamma)$ is  a line, then $\psi(\psi^{-1}(\Gamma))=\Gamma$ is a triangle by line-preservation. Otherwise $\psi^{-1}(\Gamma)$
	is the set of midpoints of a cycle in $(J,\Lambda)$ by circuit-friendliness. Since $\psi$ is line-pres it follows, as observed at the beginning of 4C5, that $\psi(\psi^{-1}(\Gamma))=\Gamma$ must be the rim-set of a wheel.

Call a triangle in a graph {\it extendible} if it is part of a {\it nondegenerate} wheel. Here comes a crisp sufficient condition for circuit-friendliness. Given a graph $G=(V,E)$ and PLS $(J,\Lambda)$ call the bijection $\psi:J\ra E$ {\it triangle-friendly} if $\psi^{-1}(\Delta)$ is a line for each extendible triangle $\Delta$ of $G$. So triangle-friendliness is a kind of dual to line-preservation. It is an exercise to verify that when $G$ is {\it wheely}, then triangle-friendly implies circuit-friendly. As to future constructions, if $J\ra E$ is line-pres (which, admittedly, is a crisper concept than its dual) then it is 'friendly' to as many triangles as there are lines. It hence seems natural to construct a fully triangle-friendly $J\ra E$ based on some preliminary line-preserving  $J\ra E$.

	{\bf 8E $M_3$-cycles.} Lemma 20 in [15] is what makes the proof of Theorem 6.2 tick. Put another way, it is because of Lemma 20 that circuit-friendliness is defined the way it is. On a technical level Lemma 20 is all about $M_3$-cycles. What are they? Each line $\ell$ of a modular lattice $L$ belongs to a {\it line-interval}. This is the interval $M_3$-sublattice of $L$ whose top is $\overline{\ell}$, and whose bottom is the meet of the three lower covers of $\overline{\ell}$ . Conversely, each interval sublattice  $[b,a]\simeq M_3$, whose top has only three lower covers in $L$, occurs this way (often hosting more than one line). The lattice in Figure 22a has five line-intervals. One of them is $[5',a_4]$ which hosts the two lines $\{1',6',8'\}$ and $\{1',7',8'\}$. Generally the top $a$ of a line-interval $[b,a]$ is called $M_3$-{\it element} (or {\it essential}) in [15,p.211]. There are subtle connections between cycles in MoPLSes (in the sense of 3A) and 'cycles of $M_3$-elements'
in the sense of [15,p.225]. Each cycle of essential elements in $L$ induces a cycle in any given MoPLS of $L$, but the converse fails. For instance the cycle $(2',3',5')$ in the MoPLS of Figure 22b doesn't yield  a cycle of the corresponding three essential elements $a_1,\ a_2,\ a_3$. 

Various types of $M_3$-elements are investigated in [15]. For instance, the existence of a MoQIMP or MoUMP implies [15, Lemma 21] that all $M_3$-elements are of  a certain type $\le (2.1)$ or $\le(3.1)$ respectively. Conversely, confining the types of occuring $M_3$-elements influences the shape of MoPLSes [15, Lemma 22]. Furthermore,
a thin lattice with a $M_3$-element of type (3.3s), such as $a_1\in L_1$ in Figure 22a, is not
tightly embeddable into a partition lattice by [15, Theorem 6].  Hence it follows from Theorem 5.4  that no graph lattice-models $L_1$. This does not trigger 'yes' for Question 3 because $L_1$ is not a BMPL, and so by Lemma 4.9(a) its MoPLS $(J_1,\Lambda_1)$ (and its other MoPLSes) cannot be line-pres rank-modeled by a graph. 

Since the proofs of [15,Lemma 20] and the supporting [15, Lemma 17]  are lengthy and technical, we ask:

{\bf Question 7:} Can the proof of [15, Lemma 20] be shortened and/or the scope of Lemma 20 be widened?

By 'widening the scope of Lemma 20' we mean that the inequality $q\le..\ (q\in X)$ in $(5.1*)$ can possibly be proven for suitable sets $X\subseteq J(L)$ {\it other} than sets of midpoints of cycles. For instance $X=\{p,q,r,s\}\subseteq J(L_2)$ for $L_2$ in Figure 21a is of this type.

	{\bf 8F Triangle configurations in modular lattices.}  As previously mentioned, Huhn  showed [6] that a modular lattice is 2-distributive iff it does not contain as an interval a complemented height 3 modular lattice which is directly indecomposable. The 2-distributivity is also equivalent [4,p.388] to the avoidance of 'triangle configurations'. Here a {\it triangle configuration} [4, p.370] in a modular lattice $L$ consists of six pairwise incomparable points (=join irreducibles) $a,b,c,p,q,x$ arranged in four lines as shown in Figure 11b. Moreover no other three of these six points must constitute a line. (It is allowed that two of the four lines determine the same line-interval.) For instance, the only MoPLS of $L=L(GF(2)^3)$ is the projective plane in Figure 11a. It is labelled in such a way that one of its many triangle configurations is again rendered by Figure 11b. 

To illustrate triangle configurations further, in Figure 10b it is shown that the BMPL $(J_3,\Lambda_3)$ of Figure 9a cannot be modeled by a binary matroid. Hence $(J_3,\Lambda_3)$ cannot occur as MoPLS in a  2-distributive modular lattice by [5, Thm 5.1]. Can it occur in a modular lattice $L$ which is  {\it not} 2-distributive?
If yes then four suitable lines somewhere in $L$ constitute a triangle configuration. (One checks that the particular MoPLS$(J_3,\Lambda_3)$ has no triangle configuration.)

{\bf 8G Locally acyclic modular lattices.} In the proof of Theorem 6.7 the construction of the bi-modeling binary  matroid $M(E)$ in [5,Thm.5.1] proceeds by induction on the height of $L$. In doing so, the MoPLS $(J,\Lambda)$ yields some 'localized' PLS which, loosely speaking, is the difference between $(J,\Lambda)$ and the  MoPLS induced by the interval sublattice $[0,a]$ where $a$ is any co-atom of $L$. It is a crucial consequence of 2-distributivity (in fact equivalent to it) that this local PLS is acyclic, even when $(J,\Lambda)$ itself is not. For instance, take the coatom $a_2$ of $L_1$ in Figure 22a. From $J(a_2)=\{3',4',5',6',7'\}$ follows that the localization at $a_2$ is the acyclic PLS  depicted in Figure 22c. (Of course this type of PLS cannot help but having 2-element lines. This is the only exception to our 3-element policy put forth at the beginning of Section 3.) In contrast take $L=L(GF(2)^3$. The seven co-atoms of $L$ match the seven lines of the MoPLS in Figure 11a. For instance the line $\{a,p,c\}$ yields the co-atom $a\vee p\vee c$.
Localizing  at the coatom given by the 'curved line' $\{p,x,q\}$ yields the PLS in Figure 11c, which is very much cyclic.

A leisurely account of [5, Thm. 5.1] is Theorem 16 in [14] where other results of [5] are recast as well. Yet [14] is a long way from being publishable. Besides tackling [5] the author ponders rendering more digestible results in [4], and perhaps some theorems of R. Wille, such as the  $FM(P,\le)$ result glanced in 2B. Collaboration on this long-term project is welcome, and may in fact be necessary to get the author going.

{\bf 8H Two loose ends.} The lattices $Fib(r)$, and also $Z(r)$, are nice examples of 2-distributive modular lattices arising in Combinatorics (Young Tableaux and the likes). Tesler [12] speaks of 'strongly modular' lattices and seems not aware that his definition is equivalent to 2-distributivity. 
As to the second loose end, recall the remarks in Subsection 1C about the super-exponentiality occuring in the Pudlak-T{\' u}ma Theorem. Hence our last question.

{\bf Question 8:} Can parts of the machinery developped in this article be adapted to handle {\it non-tight} but injective embeddings of modular lattices $L$ into moderate-size partition lattices?

A natural place to start would be to look at $L=M_4,\ M_5, \ $ and so forth. What are the smallest values of $n$ for which these lattices are embeddable into $Part(n)$?

\section*{References}
\begin{enumerate}
\item [{[1]}] U. Faigle, C. Herrmann, Projective geometry on partially ordered sets, Trans Amer. Math. Soc. 266 (1981) 319-332.
\item [{[2]}] G. Gr\"{a}tzer, Lattice Theory: Foundation, Birkh\"{a}user, 2011.
\item [{[3]}] G. Gr\"{a}tzer, E. Kiss, A construction of semimodular lattices, Order 2 (1986) 351-365.

\item [{[4]}] C. Herrmann, D. Pickering, M. Roddy, A geometric description of modular lattices, Algebra Universalis 31 (1994) 365-396.
	\item [{[5]}] C. Herrmann, M. Wild, Acyclic modular lattices and their representations, J. Algebra 136 (1991) 17-36.
	\item [{[6]}] A. Huhn, Two notes on n-distributive lattices, Colloquia Mathematics Societatis J{\'a}nos Bolyai 14 (1977) 137-147.
	\item [{[7]}] B. J{\'o}nsson and J.B. Nation, Representations of 2-distributive modular lattices of finite length, Acta Sci. Math 51 (1987) 123-128.
	\item [{[8]}] J. Mighton, A new characterization of graphic matroids, J. Comb. Theory B 98 (2008) 1253-1258.
	\item [{[9]}] J.B. Nation, Notes on Lattice Theory (unpublished).
	\item [{[10]}] J.G. Oxley, Matroid Theory, Oxford Graduate Texts in Mathematics 3, 1997.
	\item [{[11]}] P. Pudlak and J. T{\' u}ma, Every finite lattice can be embedded in a finite partition lattice, Algebra Universalis 10 (1980) 74-95.
	\item [{[12]}] G.P. Tesler, Semi-primary lattices and Tableaux algorithms, PhD, MIT 1995.
	\item [{[13]}] W.T. Tutte, Matroids and graphs, Trans Amer. Math. Soc. 90 (1959) 527-552.
		\item[{[14]}] M. Wild, Modular lattices of finite length, 28 pages, 1992, unpublished but on ResearchGate. (See also Subsection 8G.)
	\item[{[15]}] M. Wild, Cover preserving embedding of modular lattices into partition lattices, Discrete Mathematics 112 (1993) 207-244.
	\item[{[16]}] M. Wild, The minimal number of join irreducibles of a finite modular lattice, Algebra Universalis 35 (1996) 113-123.
	\item[{[17]}] R. Wille, \"{U}ber modulare Verb\"{a}nde, die von einer endlichen halbgeordneten Menge frei erzeugt werden, Math. Z. 131 (1973), 241–249. 
\end{enumerate}
\end{document}